\documentclass[11pt]{article}
\usepackage{pstricks}
\usepackage{enumerate}
\begin{document}
\newtheorem{defn0}{Definition}[section]
\newtheorem{prop0}[defn0]{Proposition}
\newtheorem{thm0}[defn0]{Theorem}
\newtheorem{lemma0}[defn0]{Lemma}
\newtheorem{coro0}[defn0]{Corollary}
\newtheorem{exa}[defn0]{Example}
\newtheorem{rem}[defn0]{Remark}

\def\rig#1{\smash{ \mathop{\longrightarrow}\limits^{#1}}}
\def\nwar#1{\nwarrow
   \rlap{$\vcenter{\hbox{$\scriptstyle#1$}}$}}
\def\near#1{\nearrow
   \rlap{$\vcenter{\hbox{$\scriptstyle#1$}}$}}
\def\sear#1{\searrow
   \rlap{$\vcenter{\hbox{$\scriptstyle#1$}}$}}
\def\swar#1{\swarrow
   \rlap{$\vcenter{\hbox{$\scriptstyle#1$}}$}}
\def\dow#1{\Big\downarrow
   \rlap{$\vcenter{\hbox{$\scriptstyle#1$}}$}}
\def\up#1{\Big\uparrow
   \rlap{$\vcenter{\hbox{$\scriptstyle#1$}}$}}
\def\lef#1{\smash{ \mathop{\longleftarrow}
    \limits^{#1}}}
\def\O{{\cal O}}
\def\L{{\cal L}}
\def\P#1{{\bf P}^#1}
\newcommand{\defref}[1]{Definition~\ref{#1}}
\newcommand{\propref}[1]{Proposition~\ref{#1}}
\newcommand{\thmref}[1]{Theorem~\ref{#1}}
\newcommand{\lemref}[1]{Lemma~\ref{#1}}
\newcommand{\corref}[1]{Corollary~\ref{#1}}
\newcommand{\exref}[1]{Example~\ref{#1}}
\newcommand{\secref}[1]{Section~\ref{#1}}
\newcommand{\eqref}[1]{(\ref{#1})}

\newcommand{\qedd}{\hfill\framebox[2mm]{\ }\medskip}
\newcommand{\codim}{\textrm{codim}}
\newcommand{\mod}{\textrm{mod}}

\author{Maria Chiara Brambilla and Giorgio Ottaviani}
\title{On partial polynomial interpolation}
\date{}
\maketitle
\abstract{The Alexander-Hirschowitz theorem says that a
general collection of $k$ double points in ${\bf P}^n$ imposes
independent conditions on homogeneous polynomials of degree $d$
with a well known list of exceptions. We generalize this theorem
to arbitrary zero-dimensional schemes contained in a general union
of double points. We work in the polynomial interpolation setting.
In this framework our main result says that the affine space of
polynomials of degree $\le d$ in $n$ variables, with assigned
values of any number of general linear combinations of first
partial derivatives, has the expected dimension if $d\neq 2$ with
only five exceptional cases. If $d=2$ the exceptional cases are
fully described.}
\medskip

Both authors are partially supported by Italian MUR and are
members of GNSAGA-INDAM.

\section{Introduction}

Let $R_{d,n}=K[x_1,\ldots ,x_n]_d$
be the vector space of polynomials of degree $\le d$ in $n$
variables over an infinite field $K$. Note that $\dim R_{d,n}={{n+d}\choose{d}}$. Let
$p_1,\ldots ,p_k\in{K}^n$ be $k$ general points and assume that
over each of these points a general affine proper subspace
$A_i\subset K^{n}\times K$ of dimension $a_i$ is given. Assume
that $a_1\geq\ldots\geq a_k$. Let $\Gamma_f\subseteq K^{n}\times
K$ be the graph of $f\in R_{d,n}$ and $T_{p_i}\Gamma_f$ be its
tangent space at the point $(p_i,f(p_i))$. Note that $\dim
T_{p_i}\Gamma_f=n$ for any $i$. Consider the conditions
\begin{equation}\label{stella}
A_i\subseteq T_{p_i}\Gamma_f, \mbox{ for }i=1,\ldots,k
\end{equation}
When $a_i=0$, the assumption \eqref{stella} means that the value
of $f$ at $p_i$ is assigned. When $a_i=n$, \eqref{stella} means
that the value of $f$ at $p_i$ and the values of all first partial
derivatives of $f$ at $p_i$ are assigned. In the intermediate
cases, \eqref{stella} means that the value of $f$ at $p_i$ and the
values of some linear combinations of first partial derivatives of
$f$ at $p_i$ are assigned.

Consider now the affine space
\begin{equation}\label{vdn}
V_{d,n}(p_1,\ldots ,p_k,A_1,\ldots, A_k)=\{f\in R_{d,n}|
A_i\subseteq T_{p_i}\Gamma_f, \quad i=1,\ldots,k\}\end{equation}
 The
polynomials in this space solve a partial polynomial interpolation
problem. The conditions in \eqref{stella} correspond to $(a_i+1)$
affine linear conditions on $R_{d,n}$. Our main result describes
the codimension of the above affine space. Since the description
is different for $d=2$ and $d\neq 2$, we divide the result in two
parts.

\begin{thm0}\label{main1}
Let  $d\neq 2$ and $\textrm{char}\,(K)=0$. 
For a general choice of points $p_i$ and subspaces
$A_i$, the affine space $V_{d,n}(p_1,\ldots ,p_k,A_1,\ldots, A_k)$
has codimension in $R_{d,n}$ equal to
$$\min\{\sum_{i=1}^k(a_i+1),\dim R_{d,n}\}$$ with the following
list of exceptions

\begin{eqnarray*}
a)& n=2,\ d=4,\ k=5,& a_i=2\textrm{ for }i=1,\ldots 5\\
b)& n=3,\ d=4,\ k=9,&  a_i=3\textrm{ for }i=1,\ldots 9\\
b')& n=3,\ d=4,\ k=9,&  a_i=3\textrm{ for }i=1,\ldots 8\textrm{ and  }a_9=2\\
c)& n=4,\ d=3,\ k=7,&  a_i=4\textrm{ for }i=1,\ldots 7\\
d)& n=4,\ d=4,\ k=14,& a_i=4\textrm{ for }i=1,\ldots 14
\end{eqnarray*}

In particular when  $\sum_{i=1}^k(a_i+1)={{n+d}\choose{d}}$ there
is a unique polynomial $f$ in $V_{d,n}(p_1,\ldots ,p_k,A_1,\ldots
A_k)$, with the above exceptions a), b'), c), d). In the
exceptional cases the space $V_{d,n}(p_1,\ldots ,p_k,A_1,\ldots
A_k)$ is empty.
\end{thm0}

The ``general choice'' assumption means that the points can be
taken in a Zariski open set (i.e.\ outside the zero locus of a
polynomial)
 and for each of these points the space $A_i$ can
be taken again in a Zariski open set. On the real numbers this
assumption means that the choices can be done outside a set of
measure zero.
  Our result is not constructive but it ensures that
in the case $\sum_{i=1}^k(a_i+1)={{n+d}\choose{d}}$
the linear system computing the interpolating polynomial with general data has a unique solution.
Hence any algorithm solving linear systems can be successfully applied.
Actually our proof shows that \thmref{main1} holds on any infinite field,
with the possible exception of finitely many values of $\textrm{char\ }K$ (see the appendix).
For finite
fields the genericity assumption is meaningless.

The case in which $a_i=n$ for all $i$ was proved by Alexander and
Hirschowitz in \cite{AH1,AH2}, see \cite{BO} for a survey. The most notable exception is the case of seven
 points with seven tangent spaces
for cubic polynomials in four variables, as in c).
 This example was known to classical algebraic geometers and it was
rediscovered in the setting of numerical analysis in \cite{Rei}.
The case of
curvilinear schemes was proved as a consequence of a more general
result by \cite{CG} on ${\bf P}^2$ and by \cite{CM} in general.

The case $d=1$ follows from elementary linear algebra. The case
$n=1$ is easy and well known: in this case the statement of
\thmref{main1} is true with the only requirement that the points
$p_i$ are distinct and the spaces $A_i$ are not vertical, that is
their projections $\pi(A_i)$ on $K^n$ satisfy $\dim
A_i=\dim\pi(A_i)$.

Assume now $d=2$. We set $a_i=-1$ for $i>k$. For any $1\le i \le
n$ we denote
$$\delta_{a_1,\ldots ,a_k}(i)=\max\{0,\sum_{j=1}^i a_j-\sum_{j=1}^i (n+1-j)\}$$

\begin{thm0}\label{main2}  Let $K$ be an infinite field. For a general choice of points $p_i$ and
  subspaces $A_i$,
the affine space $V_{2,n}(p_1,\ldots ,p_k,A_1,\ldots A_k)$ has
codimension in $R_{2,n}$ equal to $$\min\{\sum_{i=1}^k(a_i+1),\dim
R_{2,n}\}$$if and only if one of the following conditions takes
place:
\begin{enumerate}
\item
either $\delta_{a_1,\ldots ,a_k}(i)=0$ for all $1\leq i\leq n$;
\item
or $\sum_i (a_i+1)\geq {{n+2}\choose{2}}+\max\{\delta_{a_1,\ldots
  ,a_k}(i):1\leq i\leq n\}$.
\end{enumerate}

In particular when $\sum_{i=1}^k(a_i+1)={{n+2}\choose 2}$ there is
a unique polynomial $f$ in $V_{2,n}(p_1,\ldots ,p_k,A_1,\ldots
A_k)$ if and only if, for any  $1\le i\le n$, we have
$$\sum_{j=1}^i a_j\le \sum_{j=1}^i (n+1-j).$$
\end{thm0}

The first nontrivial example which explains \thmref{main2} is the
following. Consider $k=2$ and $(a_1, a_2)=(n, n)$. Then the affine
space $V_{2,n}(p_1,p_2, A_1, A_2)$  is given by quadratic
polynomials with assigned tangent spaces $A_1$, $A_2$ at two points $p_1$,
$p_2$. This space is not empty if and only if the intersection
space $A_1\cap A_2$ is not empty and its projection on $K^n$ contains
the midpoint of $p_1p_2$, which is a codimension one condition.
In order to prove this fact restrict to the line through $p_1$ and
$p_2$ and use a well known property of the tangent lines to the
parabola.
In this case
$\delta_{n,n}(i)=\left\{\begin{array}{cl}
0&i\neq 1\\
1&i=1
\end{array}\right.$
and the two conditions of \thmref{main2} are not satisfied. In
Section 3 we will explain these two conditions in graphical terms.

Let $\pi(A_i)$ be the projection of $A_i$ on $K^{n}$. For $i=1,\ldots,k$
we consider the ideal
$$I_i=\{f\in K[x_1,\ldots ,x_n]|
f(p_i)+\sum_{j=1}^n(x_j-(p_i)_j)\frac{\partial f}{\partial x_j}(p_i)=0
\mbox{ for any }x\in\pi(A_i)\}$$ Notice that we have
${m}_{p_i}^2\subseteq I_i\subseteq {m}_{p_i}$ and the ring
$K[x_1,\ldots ,x_n]/I_i$ corresponds to a
zero-dimensional scheme $\xi_i$ of length $a_i+1$,
supported at $p_i$ and contained
in the double point $p_i^2$. When
$V_{d,n}(p_1,\ldots ,p_k,A_1,\ldots, A_k)$ is not empty, its
associated vector space (that is its translate containing the
origin) consists of  the hypersurfaces of degree $d$ through
$\xi_1,\ldots ,\xi_k$.  Moreover, when this vector space has the
expected dimension, it follows that $V_{d,n}(p_1,\ldots
,p_k,A_1,\ldots, A_k)$ has the expected dimension too.

The space $K^{n}$  can be embedded in the projective space ${\bf
  P}^n$. Since the choice of points is general, we can always avoid
   the ``hyperplane at infinity''.
In order to prove the above two theorems, we will reformulate them
in the projective language of hypersurfaces of degree $d$ through
zero-dimensional schemes. More precisely we refer to Theorem
\ref{quadrics} for $d=2$, Theorem \ref{cubiche} for $d=3$ and
Theorem \ref{main} for $d\geq4$. This reformulation is convenient
mostly to rely on the wide existing literature on the subject. In
this setting Alexander and Hirschowitz proved that a general
collection of double points imposes independent conditions on the
hypersurfaces of degree $d$ (with the known exceptions) and our
result generalizes to a general zero-dimensional scheme contained
in a union of double points. It is possible to degenerate such a
scheme to a union of double points only in few cases, in such
cases of course our result is trivial from \cite{AH1}.

 Our proof of \thmref{main}, and hence of \thmref{main1}, is by
induction on $n$ and $d$.
Since it is enough to
find a particular zero-dimensional scheme which imposes
independent condition on hypersurfaces of degree $d$, we
specialize some of the points on a hyperplane, following a
technique which goes back to Terracini. We need a generalization
of the Horace method, like in \cite{AH1}, that we develop in the proof
of Theorem \ref{main}. The case of cubics,
which is the starting point of the induction, is proved by
generalizing the approach of \cite{BO}, where we restricted to a
codimension three linear subspace. This case is the crucial step
which allows to prove the Theorem \ref{main1}.
Section $4$ is devoted to this case, which requires a lot of
effort and technical details, in the setting of discrete mathematics.
Compared with
the  quick proof we gave in \cite{BO},
here we are forced to divide the proof in several cases and subcases.
While the induction argument works quite smoothly for $n, d\gg 0$, it
is painful to cover many of the initial cases. In the case $d=3$ we need
the help of a computer,
by a Montecarlo technique explained in the appendix. 

A further remark is necessary. In  \cite{AH1,BO} the
result about the independence of double points was shown to be
equivalent, through Terracini lemma, to a statement about the
dimension of higher secant varieties of the Veronese varieties, which in turn
is related to the Waring problem for
polynomials. Here the assumption that $K$ is algebraically closed
of zero characteristic is necessary to translate safely the
results, see also Theorem $6.1$ and Remark $6.3$ in \cite{IK}.
For example, on the real numbers, the closure in the
euclidean topology of the locus of secants to the twisted cubic is
a semi-algebraic set, corresponding to the cubic polynomials which
have no three distinct real roots, which is defined by the condition that
the discriminant is nonpositive. Indeed a real cubic polynomial can
be expressed as the sum of two cubes of linear polynomials (Waring
problem) if and only if it has two distinct complex conjugate
roots or a root of multiplicity three.

\section{Preliminaries}
Let $X$ be a scheme contained in a collection of double points of
$\P n$. We say that the type of $X$ is $(m_1,\ldots ,m_{n+1})$ if
$X$ contains exactly $m_i$ subschemes of a double point of length
$i$, for $i=1,\ldots ,n+1$. For example the type of $k$ double
points is $(0,\ldots,0,k)$. The degree of $X$ is $\deg X=\sum i m_i$.
A scheme of type $(m_1,\ldots ,m_{n+1})$ corresponds to a collection
of linear subspaces $L_i\subseteq \P{n}$
with $\dim L_i=i-1$ with a marked point on each $L_i$.
Algebraic families of such schemes can be defined over any field $K$
with the Zariski topology.

We recall now some notation and results from \cite{BO}.

Given a zero-dimensional subscheme $X\subseteq{\bf P}^n$,
the corresponding ideal sheaf ${\cal I}_X$
and a linear system ${\cal D}$ on ${\bf P}^n$, the Hilbert
function is defined as follows:
$$h_{{\bf P}^n}(X,{\cal D}):=\dim{\rm H}^0({\cal D})-\dim{\rm
  H}^0({\cal I}_X\otimes{\cal D}).$$
If $h_{{\bf P}^n}(X,{\cal D})=\deg X$, we say that $X$ is  ${\cal
D}$-independent, and in the case  ${\cal D}={\cal O}_{\P{n}}(d)$, we say
$d$-independent.

A zero-dimensional scheme is called curvilinear if it is contained
in a non singular curve. Notice that a curvilinear scheme contained in
a double point has length $1$ or $2$.

\begin{lemma0}[Curvilinear Lemma \cite{Ch,BO}]
\label{curvilinear}
Let $X$ be a zero-dimensional \linebreak
scheme of finite length
contained in a union of double points of ${\bf P}^n$
and $\cal D$ a linear system
on ${\bf P}^n$. Then $X$ is $\cal D$-independent if and only if
every curvilinear subscheme of $X$ is $\cal D$-independent.
\end{lemma0}

Let us denote ${\cal I}_X(d)={\cal I}_X\otimes {\cal O}(d)$ and $I_X(d)={\rm
  H}^0({\cal I}_X(d))$.
The expected dimension of the vector space $I_X(d)$  is
$\textrm{expdim} (I_X(d))= \max({{n+d}\choose{n}}-\deg X,0).$

For any scheme $X\subset{\bf P}^n$ and any hyperplane $H\subseteq \P{n}$,
the residual of $X$ with respect to $H$ is denoted by $X:H$ and it is defined by
the ideal sheaf \linebreak 
${\cal I}_{X:H}={\cal I}_{X}:{\cal I}_{H}$.
We have, for any $d$, the  well known {\em Castelnuovo sequence}
$$0\rig{}I_{(X:H)}(d-1)\rig{}I_{X}(d)\rig{}I_{(X\cap H),H}(d).$$

\begin{rem}\label{perfetto}
\textup{
If $Y\subseteq X\subseteq\P{n}$ are zero-dimensional schemes, then
\begin{itemize}
\item if $X$ is $d$-independent, then so is $Y$,
\item if $h_{{\bf P}^n}(Y,d)={{d+n}\choose{n}}$, then
$h_{{\bf P}^n}(X,d)={{d+n}\choose{n}}$.
\end{itemize}
It follows that if any zero-dimensional scheme $X\subseteq\P{n}$ with $\deg
X={{d+n}\choose{n}}$ is $d$-independent, then any
scheme contained in $X$ imposes independent conditions on
hypersurfaces of degree $d$ in ${\bf P}^n$.}
\end{rem}

\begin{rem}\label{perfetto-difettivo}
\textup{
Fix $n\geq2$ and $d\geq 3$.
Assume that if a scheme $X$ with degree ${{d+n}\choose{n}}$ does not
impose independent conditions
on hypersurfaces of degree $d$ in $\P{n}$, then it is of type
$(m_1,\ldots,m_{n+1})$ for some given $m_i$.
It follows that any subscheme of $X$ is $d$-independent.
Indeed any proper subscheme $Y$ of $X$ is also a subscheme of a
scheme $X'$ with degree ${{d+n}\choose{n}}$ and of type
$(m'_1,\ldots,m'_{n+1})\neq(m_1,\ldots,m_{n+1})$,
for some $m'_i$ and since $X'$ is $d$-independent, so is $Y$.
Moreover any scheme $Z$ containing $X$ impose independent
conditions on hypersurfaces of degree $d$ if it
contains a scheme $X''$ with degree
${{d+n}\choose{n}}$ and of type
$(m''_1,\ldots,m''_{n+1})\neq(m_1,\ldots,m_{n+1})$
for some $m''_i$. Indeed since $X''$ imposes independent conditions
on hypersurfaces of degree $d$, also $Z$ does.
}
\end{rem}

\section{Quadratic polynomials}
Assume that $X$ is a scheme of type $(m_1,\ldots,m_{n+1})$.
Let us fix an order on the irreducible components
$\xi_1,\ldots,\xi_m$ of $X$ (where $m=\sum m_i$) such that
$$\textrm{length} (\xi_1)\geq \ldots \geq \textrm{length} (\xi_{m})$$
 and for any $1\le i \le m$ let us denote by $l_i$ the length of $\xi_i$ and by $p_i$ the point
where $\xi_i$ is supported. Set $l_i=0$ for $i>m$.
For any $1\le i \le n$ let us  denote
$$\delta_X(i)=\max\{0,\sum_{j=1}^i l_j-\sum_{j=1}^i (n+2-j)\}.$$

\begin{rem}\textup{
Notice that if $\delta_X(i)$ is positive for some $1\le i \le n$, then
the scheme $\{\xi_1,\ldots,\xi_i\}$ does not impose independent
conditions on quadrics.
Indeed let be $\{e_0,\ldots,e_n\}$ a basis of $V$ and
$A$ be the symmetric matrix defining a quadric in $\mathbf{P}(V)$
passing through the scheme $\{\xi_1,\ldots,\xi_i\}$.
We may assume that $p_j=[e_{n+1-j}]$ for all $j=1,\ldots, i$.
Then, the condition $\delta_X(i)>0$ implies that the elements in the
last $i$ columns and rows of the matrix $A$ outlined below are all equal to
$0$.
Hence the scheme $\{\xi_1,\ldots,\xi_i\}$ has degree
$\sum_{j=1}^i l_j$, but imposes only $\sum_{j=1}^i
(n+2-j)={{n+2}\choose 2}-{{n-i+2}\choose 2}$ conditions on quadrics.
In fact the quadrics containing $\{\xi_1,\ldots,\xi_i\}$ are exactly all the
quadrics singular along the linear space spanned by
$\{p_1,\ldots,p_i\}$.
}
\end{rem}

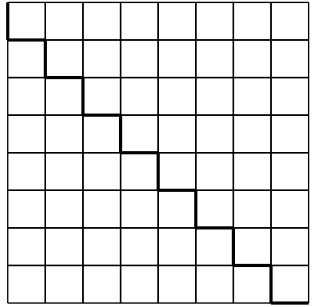
\begin{center}
\begin{pspicture}(0,0)(4,4)
\psline[linewidth=0.5pt](0,0)(4,0)\psline[linewidth=0.5pt](0,0)(0,4)
\psline[linewidth=0.5pt](4,4)(4,0)\psline[linewidth=0.5pt](4,4)(0,4)
\psline[linewidth=0.5pt](0,2)(4,2)\psline[linewidth=0.5pt](2,0)(2,4)
\psline[linewidth=0.5pt](2,4)(2,0)\psline[linewidth=0.5pt](4,2)(0,2)
\psline[linewidth=0.5pt](0,1)(4,1)\psline[linewidth=0.5pt](1,0)(1,4)
\psline[linewidth=0.5pt](1,4)(1,0)\psline[linewidth=0.5pt](4,1)(0,1)
\psline[linewidth=0.5pt](0,3)(4,3)\psline[linewidth=0.5pt](3,0)(3,4)
\psline[linewidth=0.5pt](3,4)(3,0)\psline[linewidth=0.5pt](4,3)(0,3)
\psline[linewidth=0.5pt](0,0.5)(4,0.5)\psline[linewidth=0.5pt](0.5,0)(0.5,4)
\psline[linewidth=0.5pt](0.5,4)(0.5,0)\psline[linewidth=0.5pt](4,0.5)(0,0.5)
\psline[linewidth=0.5pt](0,1.5)(4,1.5)\psline[linewidth=0.5pt](1.5,0)(1.5,4)
\psline[linewidth=0.5pt](1.5,4)(1.5,0)\psline[linewidth=0.5pt](4,1.5)(0,1.5)
\psline[linewidth=0.5pt](0,2.5)(4,2.5)\psline[linewidth=0.5pt](2.5,0)(2.5,4)
\psline[linewidth=0.5pt](2.5,4)(2.5,0)\psline[linewidth=0.5pt](4,2.5)(0,2.5)
\psline[linewidth=0.5pt](0,3.5)(4,3.5)\psline[linewidth=0.5pt](3.5,0)(3.5,4)
\psline[linewidth=0.5pt](3.5,4)(3.5,0)\psline[linewidth=0.5pt](4,3.5)(0,3.5)

\psline[linewidth=1.2pt](0,4)(0,3.5)\psline[linewidth=1.2pt](0,3.5)(0.5,3.5)
\psline[linewidth=1.2pt](0.5,3.5)(0.5,3)\psline[linewidth=1.2pt](0.5,3)(1,3)
\psline[linewidth=1.2pt](1,3)(1,2.5)\psline[linewidth=1.2pt](1,2.5)(1.5,2.5)
\psline[linewidth=1.2pt](1.5,2.5)(1.5,2)\psline[linewidth=1.2pt](1.5,2)(2,2)
\psline[linewidth=1.2pt](2,2)(2,1.5)\psline[linewidth=1.2pt](2,1.5)(2.5,1.5)
\psline[linewidth=1.2pt](2.5,1.5)(2.5,1)\psline[linewidth=1.2pt](2.5,1)(3,1)
\psline[linewidth=1.2pt](3,1)(3,0.5)\psline[linewidth=1.2pt](3,0.5)(3.5,0.5)
\psline[linewidth=1.2pt](3.5,0.5)(3.5,0)\psline[linewidth=1.2pt](3.5,0)(4,0)
\end{pspicture}
\end{center}

\medskip

The following result describes the schemes which impose
independent conditions on quadrics.
\begin{thm0}\label{quadrics}
A general zero-dimensional scheme $X\subset{\bf P }^n$ contained in a
union of double points of type $(m_1,\ldots,m_{n+1})$  imposes
independent conditions on quadrics if and only if one of the
following conditions takes place:
\begin{enumerate}
\item
either $\delta_X(i)=0$ for all $1\leq i\leq n$;
\item
or $\deg X\geq {{n+2}\choose{2}}+\max\{\delta_X(i):1\leq i\leq n\}$.
\end{enumerate}
\end{thm0}

{\it Proof.} First we prove that if $X$ does impose independent
conditions on quadrics, then either condition 1 or 2 hold. Assume
that both conditions are false and let us prove that $I_X(2)$ has
not the expected dimension $\max\{0,{{n+2}\choose2}-\deg(X)\}$.
Indeed choose $i\in\{1,\ldots, n\}$ such that $\delta_X(i)>0$ and
$\deg (X)< {{n+2}\choose{2}}+\delta_X(i)$ and consider the family
$\mathcal C$ of quadratic cones with vertex containing the linear
space ${\bf P}^{i-1}$ spanned by $p_1,\ldots,p_i$. Of course we have
$$\dim I_X(2)\geq \dim (\mathcal {C})-(\deg(X)-\sum_{j=1}^i l_j)=
{{n-k+2}\choose 2}-\deg(X)+\sum_{j=1}^i l_j=c.$$
Now we compute
$$\dim I_X(2)-\textrm{expdim} I_X(2)\geq
\min\{c,{{n-k+2}\choose 2}+\sum_{j=1}^i l_j-{{n+2}\choose2}\}=\min\{c,\delta_X(i)\}
$$
By assumption $\delta_X(i)>0$ and
$$c>{{n-k+2}\choose 2}-{{n+2}\choose 2}-\delta_X(i)+\sum_{j=1}^i l_j=
\sum_{j=1}^i l_j-\sum_{j=1}^i (n+2-j)-\delta_X(i)=0$$
Hence the dimension of $I_X(2)$ is higher than the expected dimension.

Now we want to prove that if either condition 1 or condition 2 hold, then
$X$ imposes independent conditions on quadrics.
We work by induction on $n\geq2$. If $n=2$ it is easy to check directly our claim.

Consider a scheme $X$ in $\P{n}$ which satisfies condition 1 and fix a hyperplane $H\subset\P{n}$.
We specialize all the components of $X$ on $H$ in such a way that the residual of each of
the components $\xi_1,\ldots, \xi_{n+1}$ is $1$ (if they are not zero)
and the residual of the remaining components is zero.
Then we get the Castelnuovo sequence
$$0\to I_{X:H}(1)\to I_X(2) \to I_{X\cap H}(2)$$
where $X:H$ is the residual given by at most $n+1$ simple points and $X\cap H$ is the trace in $H$.
Hence we conclude by induction once we have proved that the trace $X\cap H$ satisfies condition 1 or 2.

In fact we prove now that if $X\cap H$ does not satisfy condition 1,
then it satisfies 2.
Assume that for $X$ in $\P n$ we have
$\delta_X(i)=0$ for all $1\leq i\leq n$, while for $X\cap H$ in $H$ we have
$\delta_{X\cap H}(i)>0$ for some $1\leq i\leq n-1$. Choose the index
$1\leq k+h\leq n-1$
such that $\delta_{X\cap H}(i)$ is maximum.
This is possible only if there is some index
$1\leq k\leq k+h\leq n-1$ such that
$l_k=l_{k+1}=\ldots=l_n=l_{n+1}=\ldots=l_{n+1+h}=l$.
Since $\delta_X(k+h)=0$ and $\delta_{X\cap H}(k+h)=0$ we obtain
$$\sum_{i=1}^{k+h}(n+2-i)-h<\sum_{i=1}^{k+h} l_i\leq \sum_{i=1}^{k+h}(n+2-i),$$
from which it follows that
\begin{equation}\label{minimum-l}
l\geq (n+2-k-h).
\end{equation}
Now in order to prove that $X\cap H$ satisfies 2 we need to show that
$$\deg X\cap H \geq {{n+1}\choose{2}}+\delta_{X\cap H}(k+h).$$
Notice that
$$\deg X\cap H\geq \deg X-(n+1)\geq \sum_{i=1}^{n+1+h}l_i-(n+1),$$
hence if we prove the following inequality we are done:
$$\sum_{i=1}^{n+1+h}l_i-(n+1)\geq {{n+1}\choose{2}}+\delta_{X\cap H}(k+h)$$
i.e.
$$\sum_{i=1}^{k-1}l_i+\sum_{i=k}^{n+1+h}l-(n+1)\geq {{n+1}\choose{2}}+
\sum_{i=1}^{k-1}(l_i-1)+\sum_{i=k}^{k+h}l-\sum_{i=1}^{k+h}(n+1-i)$$
which reduces to
$$(n+1-k)l-(n+1)\geq {{n+1}\choose{2}}
-(k-1)-(k+h)(n+1)+{{k+h+1}\choose{2}}.$$
By using inequality \eqref{minimum-l} it is enough to prove, for any $n\geq 2$, any $1\leq k\leq n-1$
and any $0\leq h\leq n-1-k$,  that
\begin{equation}\label{disuguaglianza}
(n+1-k)(n+2-k-h)\geq {{n+1}\choose{2}}
+(n-k)-(k+h)(n+1)+{{k+h+1}\choose{2}}
\end{equation}
and we prove this inequality by induction on $h$.

First fix $n,k$ and choose $h=n-1-k$. In this case
\eqref{disuguaglianza} becomes
$$3(n+1-k)\geq {{n+1}\choose{2}}
+(n-k)-(n-1)(n+1)+{{n}\choose{2}}$$
which is true.
Now if we assume that \eqref{disuguaglianza} is verified for  $h=h'\leq n-1-k$,
it is easy to check it for $h=h'-1$, thus completing the proof of
\eqref{disuguaglianza}.

It remains to prove that if $X$ satisfies condition 2, then the system
of quadrics $|{\cal I}_X(2)|$
containing $X$ is empty.

If $\delta_X(i)=0$ $\forall i$ then we are in the previous case. We may assume
that there exists $i$
such that  $\delta_X(i)>0$. If the sequence $\{\delta_X(i)\}$ is nondecreasing
then
in particular $\delta_X(n)>0$. This implies easily that the quadrics containing
the first $n$ components
$\{\xi_1,\ldots ,\xi_n\}$ are singular along the hyperplane
$H=< p_1,\ldots ,p_n>$, so the only existing quadric is the double
hyperplane $H^2$.
By assumption $\deg X>\left[{{n+2}\choose 2}-1\right]+\delta_X(n)$,
hence there is
an extra condition
and   $|{\cal I}_X(2)|=\emptyset$ as we wanted.

We may assume that there exists $i<n$ such that $\delta_X(i+1)<\delta_X(i)$
and we pick the first such $i$.
In particular it follows
\begin{equation}
\label{safe}
l_{i+1}<n+1-i
\end{equation} As above, all the quadrics containing
$X_0=\{\xi_1,\ldots ,\xi_i\}$ are singular along the linear space
$L_0=<p_1,\ldots ,p_i>$ .
Let $X_1=X\setminus X_0$.
By definition $\deg X_0= {{n+2}\choose 2}-{{n+2-i}\choose 2}+\delta_X(i)$.

Let $\pi$ be the projection from $L_0$ on a linear space $L_1\simeq{\bf
P}^{n-i}$.
 By (\ref{safe}) we have $\deg X_1=\deg \pi(X_1)$.
Hence there is a bijective correspondence between $|{\cal I}_X(2)|$  and
$|{\cal I}_{\pi(X_1)}(2)|\subseteq |{\cal O}_{L_1}(2)|$.

Note that
$$\deg  \pi(X_1)-{{n+2-i}\choose 2} =\deg X-{{n+2}\choose 2}-\delta_X(i)\ge
\max_h\{\delta_X(h)\}-\delta_X(i)\ge 0$$
hence if  $\delta_{\pi(X_1)}(h)=0$ for $h=1,\ldots ,n-i$ we conclude again by
the first case.
If there exists $j$ such that  $\delta_{\pi(X_1)}(j)>0$, notice that in such
a case we have
$$\delta_{\pi(X_1)}(j)=\delta_X(j+i)-\delta_X(i),$$
hence
$$\max_p\{\delta_{\pi(X_1)}(p)\}=\max_h\delta_{\pi(X_1)}(h)-\delta_X(i). $$
We get that
$$\deg  \pi(X_1)-{{n+2-i}\choose 2}\ge  \max_p\{\delta_{\pi(X_1)}(p)\}$$

This means that $\pi(X_1)$ satisfies the assumption 2 on $L_1$
and then by (complete) induction on $n$ we get that
$|{\cal I}_{\pi(X_1)}(2)|=\emptyset$
as we wanted.
\qedd

A straightforward consequence of the previous theorem is the following
corollary.
\begin{coro0}\label{quadrics-coroll}
A general zero-dimensional scheme $X\subset{\bf P }^n$ contained in a
union of double points with $\deg X={{n+2}\choose{2}}$
imposes independent conditions on quadrics if and only if
 $\delta_X(i)=0$ for all $1\leq i\leq n$.
\end{coro0}

Theorem \ref{quadrics} provides a classification of all the types of general subschemes $X$
of a collection of double points of $\P{n}$ which do not impose
independent conditions on quadrics.
For example in $\P2$, the only case is $X$  given by two
double points.
In $\P3$ and in $\P4$ we have the following lists of subschemes which do not impose
independent conditions on quadrics.

\bigskip

\begin{center}
\begin{tabular}{|c|c|c|c|c|}
\hline
$X$ & $\deg X$ & $\max\{\delta_X(i)\}$ & $(m_1,\ldots,m_4)$ & $\dim I_X(2)$ \\
\hline
4,4,4 & 12 & 3 & $(0,0,0,3)$ & 1 \\
4,4,3 & 11 & 2 & $(0,0,1,2)$ & 1 \\
4,4,2 & 10 & 1 & $(0,1,0,2)$ & 1 \\
4,4,1,1 & 10 & 1 & $(2,0,0,2)$ & 1 \\
4,4,1 & 9 & 1 & $(1,0,0,2)$ & 2 \\
4,4 & 8 & 1 & $(0,0,0,2)$ & 3 \\
4,3,3 & 10 & 1 & $(0,0,2,1)$ & 1 \\
\hline
\end{tabular}
\\
\medskip
Table 1: List of exceptions in $\P3$
\end{center}

\bigskip

\begin{center}
\begin{tabular}{|c|c|c|c|c|}
\hline
$X$ & $\deg X$ & $\max\{\delta_X(i)\}$ & $(m_1,\ldots,m_5)$ & $\dim I_X(2)$ \\
\hline
5,5,5,5 & 20 & 6 & $(0,0,0,0,4)$ & 1 \\
\hline
5,5,5,4 & 19 & 5 & $(0,0,0,1,3)$ & 1 \\
5,5,5,3 & 18 & 4 & $(0,0,1,0,3)$ & 1 \\
5,5,5,2 & 17 & 3 & $(0,1,0,0,3)$ & 1 \\
5,5,5,1,1 & 17 & 3 & $(2,0,0,0,3)$ & 1 \\
5,5,5,1 & 16 & 3 & $(1,0,0,0,3)$ & 2 \\
5,5,5 & 15 & 3 & $(0,0,0,0,3)$ & 3 \\
\hline
5,5,4,4 & 18 & 4 & $(0,0,0,2,2)$ & 1 \\
5,5,4,3 & 17 & 3 & $(0,0,1,1,2)$ & 1 \\
5,5,4,2 & 16 & 2 & $(0,1,0,1,2)$ & 1 \\
5,5,4,1,1 & 16 & 2 & $(2,0,0,1,2)$ & 1 \\
5,5,4,1 & 15 & 2 & $(1,0,0,1,2)$ & 2 \\
5,5,4 & 14 & 2 & $(0,0,0,1,2)$ & 3 \\
5,5,3,3 & 16 & 2 & $(0,0,2,0,2)$ & 1 \\
5,5,3,2 & 15 & 1 & $(0,1,1,0,2)$ & 1 \\
5,5,3,1,1 & 15 & 1 & $(2,0,1,0,2)$ & 1 \\
5,5,3,1 & 14 & 1 & $(1,0,1,0,2)$ & 2 \\
5,5,3 & 13 & 1 & $(0,0,1,0,2)$ & 3 \\
5,5,2,2,1 & 15 & 1 & $(1,2,0,0,2)$ & 1 \\
5,5,2,2 & 14 & 1 & $(0,2,0,0,2)$ & 2 \\
5,5,2,1,1,1 & 15 & 1 & $(3,1,0,0,2)$ & 1 \\
5,5,2,1,1 & 14 & 1 & $(2,1,0,0,2)$ & 2 \\
5,5,2,1 & 13 & 1 & $(1,1,0,0,2)$ & 3 \\
5,5,2 & 12 & 1 & $(0,1,0,0,2)$ & 4 \\
5,5,1,1,1,1,1 & 15 & 1 & $(5,0,0,0,2)$ & 1 \\
5,5,1,1,1,1 & 14 & 1 & $(4,0,0,0,2)$ & 2 \\
5,5,1,1,1 & 13 & 1 & $(3,0,0,0,2)$ & 3 \\
5,5,1,1 & 12 & 1 & $(2,0,0,0,2)$ & 4 \\
5,5,1 & 11 & 1 & $(1,0,0,0,2)$ & 5 \\
5,5 & 10 & 1 & $(0,0,0,0,2)$ & 6 \\
\hline
5,4,4,2 & 15 & 1 & $(0,1,0,2,1)$ & 1 \\
5,4,4,1,1 & 15 & 1 & $(2,0,0,2,1)$ & 1 \\
5,4,4,1 & 14 & 1 & $(1,0,0,2,1)$ & 2 \\
5,4,4 & 13 & 1 & $(0,0,0,2,1)$ & 3 \\
\hline
4,4,4,4 & 16 & 2 & $(0,0,0,4,0)$ & 1 \\
4,4,4,3 & 15 & 1 & $(0,0,1,3,0)$ & 1 \\
\hline
\end{tabular}
\\
\medskip
Table 2: List of exceptions in $\P4$
\end{center}

\section{Cubic polynomials}

In this section we generalize the approach of \cite[Section 3]{BO}
to our setting and we prove the following result.

\begin{thm0}\label{cubiche}
A general zero-dimensional scheme $X\subset{\bf P }^n$ contained in a  union of double points
imposes independent conditions on cubics
with the only exception of $n=4$ and $X$ given by $7$ double points.
\end{thm0}

First we give the proof of the previous theorem in cases $n=2,3,4$.

\begin{lemma0}\label{casi-iniziali}
Let be $n=2,3$ or $4$. Then a general zero-dimensional scheme
$X\subset{\bf P }^n$ contained in a union of double points
imposes independent conditions on cubics with the only
exception of $n=4$ and $X$ given by $7$ double points.
\end{lemma0}
{\it Proof.}
By Remark \ref{perfetto} it is enough to prove the statement for $X$
with degree ${{n+3}\choose 3}$.

Let $n=2$ and $X$ a subscheme of a collection of double points
with $\deg X=10$. Fix a line $H$ in $\P2$ and
consider the Castelnuovo exact sequence
$$ 0 \to I_{{X}:H}(2) \to I_X(3) \to I_{X\cap H}(3). $$
It is easy to prove that it is always possible to specialize some
components of $X$ on $H$ so that $\deg(X\cap H)=4$ and that the residual
$X:H$ does not contain two double points. The last condition ensures that
$\delta_{X:H}(i)=0$ for $i=1,2$. Hence we conclude by \corref{quadrics-coroll}.

In the case $n=3$, the scheme $X$ has degree $20$.
Since there are no cubic surfaces with five singular points we can assume that
$X$ contains at most three double points. Indeed if $X$ contains $4$ double points
we can degenerate it to a collection of $5$ double points.
€We fix a plane $H$ in $\P3$ and
we want to specialize some components of $X$ on $H$ so that
$\deg(X\cap H)=10$ and that the residual $X:H$ imposes independent conditions on quadrics.
By looking at Table 1, since $\deg(X:H)=10$, it is enough to require that $X:H$
is not of the form $(0,1,0,2),(2,0,0,2)$ or $(0,0,2,1)$.
It is easy to check that this is always possible: indeed specialize on $H$ the
components of $X$ starting from the ones with higher length and
keeping the residual
as minimal as possible until the degree of the trace is $9$ or
$10$. If the degree of the
trace is $9$ and there is in $X$ a component with length $1$ or $2$
we can obviously complete
the specialization.
The only special case is given by $X$ of type $(0,0,3,2)$ and in this
case we specialize on $H$
the two double points and two components of length $3$ so that each of
them has residual $1$.

If $n=4$ the case of $7$ double points is exceptional. Assume that $X$
has degree $35$ and
contains at most $6$ double points.
We fix a hyperplane $H$ of $\P4$ and
we want to specialize some components of $X$ on $H$ so that
$\deg(X\cap H)=20$ and that the residual $X:H$
imposes independent conditions on quadrics.
By looking at Table 2, it is enough to require that $X:H$
does not contain two double points, does not contain one double point
and two components of length $4$ and it is not of the form
$(0,0,1,3,0)$.
It is possible to satisfy this conditions by specializing the
components of $X$ in the following way:
we specialize the components of $X$ on $H$ starting from the ones with
higher length and keeping the residual
as minimal as possible until the degree of the trace is maximal and
does not exceed $20$.
Then we add some components allowing them to have residual $1$ in
order to reach the degree $20$.
It is possible to check that this construction works,
except for the case $(0,0,5,0,4)$ where we have to specialize on $H$
all the double
points and $2$ of the components with length $3$ so that both have
residual $1$.
It is easy also to check that following the construction above
the residual has always the desired form, except for $X$ of the form
$(0,0,1,8,0)$, where the above rule gives a residual of type $(0,0,1,3,0)$.
In this case we make a specialization ad hoc: for example we can put
on $H$ six components of length $4$
and the unique component of length $3$ in such a way that all them
have residual $1$ and we
obtain a residual of type $(7,0,0,2,0)$ which is admissible.

Now we have to check the schemes either contained in $7$ double
points or containing $7$ double points. But this follows immediately
by Remark \ref{perfetto-difettivo}.
\qedd

We want to restrict a zero dimensional scheme $X$ of $\P n$ 
to a given subvariety $L$.
We could define the residual $X:L$ as a subscheme
of the blow-up of $\P n$ along $L$ (\cite{AH3}), but we prefer to consider
$\deg (X:L)$ just as an integer associated to $X$ and $L$.
More precisely given a subvariety $L\subset{\bf P}^n$, 
we denote $\deg (X:L)= \deg X-\deg (X\cap L)$.
In particular we will use this notion in the following cases:
$$\deg (X:L), \quad\deg (X:(L\cup M)),\quad\deg (X:(L\cup M\cup N))$$
where  $L, M, N \subset{\bf P}^n$ are three general subspaces of codimension three.
We also recall that
$$\deg (X\cap (L\cup M))=\deg (X\cap L) +\deg (X\cap M)-\deg (X\cap (L\cap M))$$
and
$$\deg X\cap (L\cup M\cup N)=\deg( X\cap L) +\deg (X\cap M)+\deg (X\cap N)-\deg( X\cap L\cap M)$$
$$-\deg (X\cap L\cap N)-\deg (X\cap M\cap N)+\deg (X\cap L\cap M\cap N).$$

The proof of Theorem \ref{cubiche}
relies on a preliminary description, 
which is inspired to the approach of \cite{BO}. 
More precisely the proof is structured as follows:
\begin{itemize} 
\item[-]
in Proposition \ref{proprep} below we generalize \cite[Proposition 5.2]{BO}, 
\item[-]
in \propref{proprep2} and \propref{proprep2bis} we generalize \cite[Proposition 5.3]{BO}, 
\item[-]
the analogue of \cite[Proposition 5.4]{BO} is contained in  
\lemref{primo}, \lemref{secondo}, \lemref{terzo} and \propref{proprep3}.
\end{itemize}

\begin{prop0}\label{proprep}
Let $n\ge 8$ and let $L, M, N \subset\P n$ be general subspaces of
codimension $3$. Let $X=X_L\cup X_M\cup X_N$ be a general scheme contained in
a union of double points, where $X_L$ (resp. $X_M$, $X_N$) is
supported on $ L$ (resp. $M$, $N$),
 such that the triple
$\left(\deg (X_L: L), \deg( X_M:M), \deg( X_N:N)\right)$
is one of the following
\begin{enumerate}[(i)]
\item $(6,9,12)$
\item $(3,12,12)$
\item $(0,12,15)$
\item $(6,6,15)$
\item $(0,9,18)$
\end{enumerate}
then there are no
cubic hypersurfaces in $\P n$ which contain $ L\cup M\cup{ N}$ and
which contain $X$.
\end{prop0}

{\it Proof.}
For $n=8$ it is an explicit computation, which can be easily performed
with the help of a computer (see the appendix).

For $n\ge 9$ the statement follows by induction on $n$.
Indeed if $n\ge 8$ it is easy to check that there are no quadrics
containing ${L}\cup{M}\cup{N}$.
Then given a general hyperplane
$H\subset \P n$ the Castelnuovo sequence induces the isomorphism
$$0\rig{}I_{{L}\cup{ M}\cup{ N},{\P n}}(3)\rig{}I_{\left({ L}\cup{ M}\cup{ N}\right)\cap H,H}(3)\rig{} 0$$
hence specializing the support of $X$ on the hyperplane $H$, since the space \linebreak
$I_{{ L}\cup{ M}\cup{ N},{\P n}}(2)$ is empty, we get
$$0 \rig{}I_{X\cup { L}\cup{ M}\cup{ N},{\P n}}(3)\rig{}I_{\left(X\cup
    { L}\cup{ M}\cup{ N}\right)\cap H,H}(3)$$
then our statement immediately follows by induction.
\qedd

\begin{rem}\label{assumptions}
\textup{
It seems likely that the previous proposition holds with much more general
assumption. Anyway the general assumption
$\deg (X_L: L)+ \deg (X_M:M) +\deg (X_N:N)=27$
is too weak, indeed the triple $(0,6,21)$ cannot be added to the list of the \propref{proprep}.
Indeed there are two independent cubic hypersurfaces in ${\bf P}^8$,
containing $L$, $M$, $N$, two general double points on $M$
and seven general double points on $N$,
as it can be easily checked  with the help of a computer (see the appendix).
 Quite surprisingly, the triple
$(0,0,27)$ could be added to the list of the \propref{proprep}, and we think
that this phenomenon has to be better understood.
In \propref{proprep} we have chosen exactly  the assumptions
that we will need in the following propositions, in order to minimize
the number of the initial checks.
}\end{rem}

For the specialization technique we need the following two easy remarks.

\begin{rem}\label{spec3}
\textup{
Let $L$, $N$ be two codimension three subspaces of $\P n$, for
$n\ge 5$.
Let $\xi$ be a general scheme contained in a double point $p^2$ supported on $L$ such that
$\deg (\xi:L)=a$, $0\le a\le 3$.  Then there is a specialization $\eta$ of $\xi$
such that the support of $\eta$ is on $L\cap N$, $\deg (\eta:L)=a$
and $\deg(\left(\eta\cap N\right):(L\cap N))=a$.
}\end{rem}

\begin{rem}\label{spec33}
\textup{
Let $L$ be a codimension three subspaces of $\P n$.
Let $X$ be a scheme contained in a double point $p^2$.
\begin{enumerate}[i)]
\item
If $\deg X=n+1$ then there is a specialization $Y$ of $X$ which is supported at $q\in L$ such that
$\deg (Y:L)=3$.
\item
If $\deg X=n$ then there are two possible specializations $Y$ of $X$ which are supported at $q\in L$
such that $\deg (Y:L)=3$ or $2$.
\item
If $\deg X=n-1$ then there are three possible specializations $Y$ of $X$ which are supported at $q\in L$
such that $\deg( Y:L)=3$, $2$ or $1$.
\item
If $\deg X\le n-2$ then there are four possible specializations $Y$ of $X$ which are supported at $q\in L$
such that $\deg (Y:L)=3$, $2$, $1$ or $0$.
\end{enumerate}
}
\end{rem}

\begin{prop0}\label{proprep2}
Let $n\ge 5$ and let ${L, M}\subset\P n$ be subspaces of codimension three.
Let $X=X_L\cup X_M\cup X_O$ be a scheme contained in a union of double points such that
$X_L$ (resp. $X_M$) is supported on $L$ (resp. $M$) and it is general among the schemes supported on $L$
(resp. $M$) and $X_O$ is general.
Assume that the following further conditions hold:
\begin{eqnarray*}
&\deg (X_L:L)+\deg( X_M:M)+\deg X_O=9(n-1),\\
&n-2\le \deg( X_L:L)\le \deg (X_M:M)\le 4n-6,\\
&3n+3\le\deg X_O\le 3n+6.
\end{eqnarray*}
Then there are no cubic hypersurfaces in $\P n$
which contain ${L}\cup{M}$ and which contain $X$.
\end{prop0}

{\it Proof.}
For $n=5,6,7$ it is an explicit computation (see the appendix). 

For $n\ge 8$, the statement follows by induction from $n-3$ to $n$.
Indeed given a third general codimension three subspace ${N}$, we get the exact sequence
$$0\rig{}I_{{ L}\cup{ M}\cup{ N},{\P n}}(3)\rig{}
I_{{ L}\cup{M},{\P n}}(3)\rig{}I_{\left({ L}\cup{ M}\right)\cap { N},{ N}}(3)\rig{}0$$
where the dimensions of the three spaces in the sequence are respectively $27$,  $9(n-1)$ and $9(n-4)$.

We will specialize now some components of $X_L$ on $L\cap N$
and some components of $X_M$ on $M\cap N$.
We denote by
$X_L'$ the union of the components of $X_L$ supported on $L\setminus N$ and
by $X_L''$ the union of the components of $X_L$ supported on $L\cap N$.
Since $n\ge 5$ we may assume also that $\deg (X_L'':(L\cup N))=0$.
Analogously let $X_M'$ and $X_M''$ denote the corresponding subschemes of $X_M$.
Now we describe more explicitly the specialization.

From the assumption $$3n+3\le \deg X_O\le 3n+6$$
it follows that in particular $X$ has at least three irreducible components
and so we may specialize all the components of $X_O$ on $N$
in such a way that $\deg (X_O:N)=9$.

Notice that the degree of the trace $X_O\cap N=\deg X_O-9$ satisfies the same inductive hypothesis
$$3(n-3)+3\le \deg (X_O\cap N)\le 3(n-3)+6$$
and we have
$$6n-15\le \deg (X_L:L)+\deg (X_M:M)\le 6n-12$$

If $\deg (X_M:M)\le 3n$,
 by using that 
$$\deg(X_L:L)\le \frac{1}{2} \left(\deg (X_L:L)+\deg (X_M:M)\right)\le \deg(X_M:M)$$
we get
$$3n-7\le \deg( X_M:M)\le 3n$$ 
$$3n-15\le \deg (X_L:L)\le 3n-6$$ 
then we can specialize $X_M$ and $X_L$
in such a way that
$\deg (X_M':M)=12$ and $\deg (X_L':L)=6$, indeed the conditions
$$n-5\le \deg (X_M:M)-12\le 4n-18$$
$$n-5\le \deg (X_L:L)-6\le 4n-18$$
are true for $n\ge 8$ and guarantee that the inductive assumptions are true
on the trace.

Now if $\deg (X_M:M)\ge 3n+1$, we have
$$3n+1\le \deg (X_M:M)\le 4n-6$$
$$2n-9\le \deg (X_L:L)\le 3n-13$$
and we can specialize in such a way that $\deg (X_L':L)=0$ and  $\deg (X_M':M)=18$.
Indeed we have, for $n\ge 6$
$$n-5\le \deg (X_M:M)-18\le 4n-18$$
$$n-5\le \deg (X_L:L)\le 4n-18$$

In any of the previous cases, the residual satisfies the assumptions of \propref{proprep},
while the trace $\left(X\cup{ L}\cup{ M}\right)\cap { N}$ satisfies the inductive assumptions on
$N=\P{{n-3}}$.
In conclusion by using the sequence
$$0\rig{}I_{X\cup {L}\cup{ M}\cup{ N},{\P n}}(3)\rig{}
I_{X\cup { L}\cup{M},{\P n}}(3)\rig{}I_{\left(X\cup{ L}\cup{ M}\right)\cap { N},{ N}}(3)$$
we complete the proof.
\qedd

The following proposition is analogous to the previous one, with a different
assumption on $\deg X_O$. In this case we need an extra assumption
on $X_L$ and $X_M$, namely that in one of them 
there are enough irreducible components with residual different from $2$.
The reason for this choice is that it makes possible to find a suitable specialization
with  residual $3$, $9$ or $15$, by the Remark
\ref{spec3} (if all the components have residual $2$, this should not be possible).

From now on we denote by $X_L^i$ (resp.\ $X_M^i$) for $i=1,2,3$  the union of the irreducible components $\xi$ of
$X_L$ (resp.\ $X_M$) such that $\deg(\xi:L)=i$ (resp.\ $\deg(\xi:M)=i$).

\begin{prop0}\label{proprep2bis}
Let $n\ge 5$ and let ${L, M}\subset\P n$ be subspaces of codimension three.
Let $X=X_L\cup X_M\cup X_O$ be a scheme contained in a union of double points such that
$X_L$ (resp. $X_M$) is supported on $L$ (resp. $M$) and it is general among the schemes supported on $L$
(resp. $M$) and $X_O$ is general.  
Assume that either the number of the irreducible components of
$X_L^1\cup X_L^3$,
or that the number of the irreducible components of
$X_M^1\cup X_M^3$ is at least $\frac{n-2}{3}$.
Assume that the following further conditions hold:
\begin{eqnarray*}
&\deg (X_L:L)+\deg( X_M:M)+\deg X_O=9(n-1),\\
&n-2\le \deg( X_L:L)\le \deg (X_M:M)\le 4n-6,\\
&3n+7\le\deg X_O\le 5n+2.
\end{eqnarray*}
Then there are no cubic hypersurfaces in $\P n$
which contain ${L}\cup{M}$ and which contain $X$.
\end{prop0}
  
{\it Proof.}
For $n=5,6,7$ it is an explicit computation (see the appendix), and the thesis is true even without the assumption
on $X_L^1\cup X_L^3$.

For $n\ge8$ the statement follows by induction from $n-3$ to $n$, by using possibly also \propref{proprep2}.
As in the previous proof, given a third general codimension three subspace ${N}$, we get the exact sequence 
$$0\rig{}I_{{ L}\cup{ M}\cup{ N},{\P n}}(3)\rig{}
I_{{ L}\cup{M},{\P n}}(3)\rig{}I_{\left({ L}\cup{ M}\right)\cap { N},{ N}}(3)\rig{}0$$

We will specialize now some components of $X_L$ on $L\cap N$
and some components of $X_M$ on $M\cap N$.
We use the same notations as in the previous proof, and
we describe more precisely the specialization in the following two cases.

\begin{enumerate}

\item Assume first that
$$3n+7\le \deg X_O\le 4n+7$$
In particular $X$ has at least four irreducible components
and we may specialize all the components of $X_O$ on $N$
in such a way that $$\deg ((X_O\cap N):N) =12$$
and so we have
$$5n-16 \le \deg (X_L:L)+\deg (X_M:M)\le 6n-16$$

In particular it follows
$$\frac{5n}{2}-8\le \deg (X_M:M)\le 4n-6$$
$$n-2\le \deg (X_L:L)\le 3n-8$$

We divide into two subcases.

In the first one we assume that the number of the irreducible components of
$X_L^1\cup X_L^3$ is at least $\frac{n-2}{3}$.
In this case we can specialize $X_M$ and $X_L$
in such a way that
$\deg (X_M':M)=12$ and $\deg (X_L':L)=3$.
Moreover there exists a specialization such that $X_L''$ has at least
$\frac{n-5}{3}=\frac{n-2}{3}-1$ components with residual $1$ or $3$. 
Indeed in $X_L'$ we keep at most one of these
components, and if we are forced to keep three components of length one,
it means that there are no components of length $2$ in $X_L$, which
implies our claim.

Notice that the conditions
$$n-5\le \deg (X_M:M)-12\le 4n-18$$
$$n-5\le \deg (X_L:L)-3\le 4n-18$$
are true for $n\ge 10$. 
They are also true for $n\ge 8$ as soon as 
$\deg (X_M:M)\ge n+7$, so we need only to check the cases 
$8\le n \le 9$ and  $\deg (X_M:M)\le n+6$, which implies
$\deg (X_L:L)\ge 4n-22$.

In this case we specialize $X_M$ and $X_L$
in such a way that $\deg (X_M':M)=6$, $\deg (X_L':L)=9$ and 
$X_L''$ has at least
$\frac{n-5}{3}=\frac{n-2}{3}-1$ components with residual $1$ or $3$. 
The conditions
$$n-5\le \deg (X_M:M)-6\le 4n-18$$
$$n-5\le \deg (X_L:L)-9\le 4n-18$$
are true if $n=9$ or if $n=8$ and 
$\deg (X_L:L)\ge n+4$.

So the remaining cases to be considered are when
$n=8$, $\deg (X_M:M)\le n+6=14$,  and $\deg (X_L:L)\le n+3=11$, 
that is when
the triple $$(\deg (X_L:L), \deg (X_M:M), \deg X_O)$$ is one of the following:
$( 10, 14, 39)$, $( 11, 13, 39)$, $( 11, 14, 38)$,
which have been checked with random choices (see the appendix) with a computer. 

In the second subcase, 
we know that the number of the irreducible components of
$X_M^1\cup X_M^3$ is at least $\frac{n-2}{3}$.
Then we can specialize $X_M$ and $X_L$
in such a way that
$\deg (X_M':M)=9$ and $\deg (X_L':L)=6$.
As above it is easy to check that 
there exists a specialization such that $X_M''$ has at least
$\frac{n-5}{3}=\frac{n-2}{3}-1$ components with residual $1$ or $3$. 

Notice that the conditions
$$n-5\le \deg (X_M:M)-9\le 4n-18$$
$$n-5\le \deg (X_L:L)-6\le 4n-18$$
are true for $n\ge 8$ as soon as one of the following conditions is satisfied 
\begin{itemize}
\item [(a)]
$\deg (X_M:M)\le 4n-17$ , which implies  
$\deg (X_L:L)\ge n+1$.
\item [(b)] $n=8$, $\deg (X_L:L)\ge n+1=9$, which implies  $\deg (X_M:M)\le 5n-17=23$
\end{itemize}
Assume then  that (a) and (b) are not satisfied.

We have 
$4n-16\le \deg (X_M:M)\le 4n-6$  
and we specialize $X_M$ and $X_L$
in such a way that $\deg (X_M':M)=15$ and $\deg (X_L':L)=0$.
The conditions
$$n-5\le \deg (X_M:M)-15\le 4n-18$$
$$n-5\le \deg (X_L:L)\le 4n-18$$
are true for $n\ge 9$ or if $n=8$ and $\deg (X_M:M)\ge n+10$.

So the remaining cases to be considered are when $n=8$, $4n-16=16\le \deg (X_M:M)\le n+9=17$
and (by case (b)) $\deg (X_L:L)\le 8$.   
The only remaining case are
$$(\deg (X_L:L), \deg (X_M:M), \deg X_O)=(7,17,39), (8,16,39), (8,17,38) $$
which we have checked with a computer.

\item Assume now that
$$4n+8\le \deg X_O\le 5n+2$$
which implies
$$4n-11\le \deg (X_L:L)+\deg (X_M:M)\le 5n-17$$
In particular $X$ has at least five irreducible components
and we may specialize all the components of $X_O$ on $N$
in such a way that $\deg ((X_O\cap N):N)=15$.

In this case we have
$$2n-5\le \deg (X_M:M)\le 4n-6$$
$$n-2\le \deg (X_L:L)\le \frac{5n-17}{2}$$
and we can specialize $X_M$ and $X_L$
in such a way that
$\deg (X_M':M)=12$ and $\deg (X_L':L)=0$.
Notice that the conditions
$$n-5\le \deg (X_M:M)-12\le 4n-18$$
$$n-5\le \deg (X_L:L)\le 4n-18$$
are true for $n\ge 12$ and also for $n\ge 8$ as soon as
$ \deg (X_M:M)\ge n+7$.

Assume now that $8\le n\le 11$ and  $\deg (X_M:M)\le n+6$, which implies 
$\deg (X_L:L)\ge 3n-17$.

In this case we specialize $X_M$ and $X_L$
in such a way that $\deg (X_M':M)=6$ and $\deg (X_L':L)=6$.
The conditions
$$n-5\le \deg (X_M:M)-6\le 4n-18$$
$$n-5\le \deg (X_L:L)-6\le 4n-18$$
are true for $n\ge 9$ and also for $n=8$ if
$ \deg (X_L:L)\ge n+1$.

The only remaining cases to be considered are then

$n=8$,   $7\le \deg (X_L:L)\le 8$,
and $ \deg (X_M:M)\le n+6= 14$
that is when the triple 
$$(\deg (X_L:L), \deg (X_M:M), \deg X_O)$$ is one of the following:
$( 7, 14, 42)$ ,
$( 8, 13, 42)$ ,
$( 8, 14, 41)$ which we have checked with a computer.

\end{enumerate}

In conclusion in any the previous cases we conclude by using the sequence
$$0\rig{}I_{X\cup {L}\cup{ M}\cup{ N},{\P n}}(3)\rig{}
I_{X\cup { L}\cup{M},{\P n}}(3)\rig{}I_{\left(X\cup{ L}\cup{ M}\right)\cap { N},{ N}}(3)$$
since the trace $\left(X\cup{ L}\cup{ M}\right)\cap { N}$ satisfies the inductive assumptions on
$N=\P{{n-3}}$ and the residual satisfies the hypotheses of \propref{proprep}.
\qedd

Let $X_O\subset\P n$ be a scheme, contained in a union of double points,
of degree $(n+1)^2+\alpha$ with $0\le\alpha\le n-1$
and $M$ be a subspace of codimension three.
Assume that $n\ge 8$ and that $X_O$ contains at most one component of degree $\le 3$.
Let $h_i$ be the number of components of $X_O$ of degree $i$ for $i=4,\ldots, n+1$
and let $h$ ($0\le h\le 3$) be the degree of the component of $X_O$  of degree $\le 3$.
Note that $\sum_{i=4}^{n+1}ih_i+h = (n+1)^2+\alpha$.
Let us choose an order on the irreducible components of $X_O$ in such a way the length
of any component is non increasing.

We consider one of the following two specializations   $X_O=X_O'\cup X_M$
where $X_M$ is supported on $M$ and $X_O'$ is supported outside $M$:

$(a)$ we choose as $X_O'$ the union of the irreducible components of $X_O$,
starting from the ones with maximal length, in such a way that
$\deg X_O'= 3(n+1)+\beta \ge 3(n+1)+\alpha$ and it is minimal.
By construction $0\le\beta-\alpha\le n$.
Let $a_i$ be the number of components of $X_M=X_O\setminus X_O'$
of degree $i$ for $i=4,\ldots, n+1$.
Then
$$\sum_{i=4}^{n+1}i a_i + h =\deg(X_M)= (n+1)(n-2)+\alpha-\beta$$

$(\widehat{a})$ we choose as $X_O'$ the union of the irreducible components of $X_O$,
starting from the ones with maximal length, in such a way that
$\deg X_O'= 3(n+1)+\widehat{\beta} \ge 3(n+1)$ and it is minimal.
By construction $0\le\widehat{\beta} \le n-1$.
Let $\widehat{a}_i$ be the number of components of $X_M=X_O\setminus X_O'$
of degree $i$ for $i=4,\ldots, n+1$.
Then
$$\sum_{i=4}^{n+1}i \widehat{a}_i + h =\deg(X_M)= (n+1)(n-2)+\alpha-\widehat{\beta}$$

In both the specializations let us denote:
$\gamma=\deg (X_M\cap M)-(n-2)^2$
and note that we have some freedom to specialize $X_M$ on $M$,
according to Remark \ref{spec33}.
If we have a specialization
with $\deg (X_M\cap M)=p$ and another specialization with
$\deg (X_M\cap M)=q$
then for any value between $p$ and $q$ there is a suitable specialization
such that $\deg (X_M\cap M)$ attains that value. We will use often this technique by evaluating
the maximum (resp. the minimum) possible value of $\deg (X_M\cap M)$ under a specialization.

\begin{lemma0} \label{primo}
If in the specialization $(a)$ we have
$$a_n+2a_{n-1}+3\sum_{i=4}^{n-2}a_i\le 1$$
then we have $a_{n+1}\neq 0$ and there exists a
specialization of type $(a)$ such that
$\gamma=\alpha\le n-4$.
\end{lemma0}
{\it Proof.}
From the assumptions it follows that $a_{j}=0$ for any $j=4,\ldots,n-1$ and
$a_n=i$ with $0\le i\le 1$.
Then $X_M$ consists of points of maximal length $n+1$ with at most one component of of length $h$ and at most
one component of length $n$. Hence $X_O'$ consists only of double points,
hence we get $\beta=0 \textrm{\ (mod\ } n+1)$,
hence we have $a_{n+1}=\frac{(n+1)(n-2)+\alpha-\beta-h-in}{n+1}$,
which is an integer, so that $\frac{\alpha-h+i}{n+1}$
is an integer, so that $\alpha=h-i\le n-4$.

It follows that  $a_{n+1}=n-2-i$,
hence
the maximum degree of $X_M\cap M$ is
$(n-2)^2+h$, the minimum degree is $(n-2-i)(n-2)+i(n-3)+(h-1)=(n-2)^2+(h-i-1)$,
and we can choose $\gamma=h-i=\alpha$.
\qedd

\begin{lemma0} \label{secondo}
If in the specialization $(a)$ we have
$$3a_{n+1}+2a_{n}+a_{n-1}\ge 3n-7+\alpha-\beta$$
then there exists a specialization of type $(\widehat{a})$
such that either $\gamma=\alpha\le n-4$ or $\gamma=\alpha-3\le n-4$.
\end{lemma0}

{\it Proof.}
Assume first $a_{n+1}=0$. Since $\alpha-\beta\ge-n$, from the assumption it follows
$$2a_n+a_{n-1}\ge 2n-7$$
Notice also that
$$a_n+a_{n-1}\le \frac{(n+1)(n-2)+\alpha-\beta}{n}\le n-2+\frac{n-2}{n}$$
hence
$$a_n+a_{n-1}\le n-2.$$
These two conditions imply that we have only the following possibilities:
$$(a_n,a_{n-1})\in\{(n-2,0)(n-3,0),(n-4,1),(n-3,1),(n-4,2),(n-5,3)\}$$

In all these cases, by performing the specialization of type $(\widehat{a})$, we have
$n-3\le \widehat{\beta}\le n-1$ or $\widehat{\beta}=0$.
Moreover it is easy to check that
$\widehat{a}_{n}=a_{n}$ if $\alpha\le\widehat{\beta}$,
$\widehat{a}_{n}=a_{n}+1$ if $\alpha>\widehat{\beta}$,
 and $\widehat{a}_{j}=a_{j}$ for any $j\le n-1$.
In any case the difference $\delta$ between the maximum
degree  of the trace $X_M\cap M$ and the minimum degree satisfies
$$\delta \ge \widehat{a}_n+2\widehat{a}_{n-1}+3\sum_{i=4}^{n-2} \widehat{a}_i+\max\{h-1,0\}.$$
We have
$\deg(X_M)=\sum_{i=4}^{n}i \widehat{a}_i + h = (n+1)(n-2)+\alpha-\widehat{\beta}$
and so
$$\sum_{i=4}^{n-2}i \widehat{a}_i +h\ge (n+1)(n-2)-\widehat{\beta}-n\widehat{a}_n-(n-1)\widehat{a}_{n-1}.$$

In the first two cases, where $(a_n,a_{n-1})=(a,0)$ and $n-3\le a\le n-2$,
we assume first  $n-3\le \widehat{\beta}$, then
the maximal degree of the trace $X_M\cap M$ is
$$(n-2)^2+\alpha+1\le (n+1)(n-2)+\alpha-\widehat{\beta}-2\widehat{a}_n\le (n-2)^2+\alpha+3$$
since $\widehat{a}_n\ge n-3$, moreover $\delta\ge n-2\ge 6$ and so we have that
either $\gamma=\alpha$, or $\gamma=\alpha-3$ work.
It remains the case  $\widehat{\beta}=0$ where we get
that in $X_O'$ we have three points of length $n+1$, then either
$\beta=0$ and $\alpha=0$, or $\beta=n$ and $\alpha>0$.
By substituting in the hypothesis of our lemma
the values $(a_{n+1},a_n,a_{n-1})=(0,a,0)$ we get $\beta=n$ and
$0<\alpha\le 3$. In this case
 the maximal degree ${\mathcal M}$
of the trace $X_M\cap M$
satisfies
$$(n+1)(n-2)+\alpha+(n-4) \le {\mathcal M} \le (n-2)^2+\alpha+(n-2)$$
and, since $\delta\ge n-2$, the choice $\gamma=\alpha$ works.

Now consider the case $(a_n,a_{n-1})=(a,1)$, where $n-4\le a \le n-3$.
Assume first  $n-3\le \widehat{\beta}$, then the maximal degree of the trace $X_M\cap M$ is
$$(n-2)^2+\alpha\le (n+1)(n-2)+\alpha-\widehat{\beta}-2\widehat{a}_n-1\le (n-2)^2+\alpha+4$$
since $n-4\le\widehat{a}_n\le n-2$,
moreover  $\delta\ge n-1\ge 7$
so that either $\gamma=\alpha$, or $\gamma=\alpha-3$ work.
It remains the case  $\widehat{\beta}=0$, where we have
either $\beta=0$ and $\alpha=0$, or $\beta=n$ and $\alpha>0$.
By substituting in the hypothesis of our lemma
the values $(a_{n+1},a_n,a_{n-1})=(0,a,1)$, for $n-4\le a\le n-3$, we get
$\beta=n$ and $0<\alpha\le 2$. Then we have $\widehat{a}_n=n-2$ and so
 the maximal degree of the trace $X_M\cap M$ is
$$(n+1)(n-2)+\alpha-2(n-2)-1=(n-2)^2+\alpha+(n-3)$$
and since the difference $\delta\ge n-1$, the choice $\gamma=\alpha$ works.

In the case $(a_n,a_{n-1})=(n-4,2)$, if $n-3\le \widehat{\beta}$, then
 the maximal degree of the trace $X_M\cap M$ is
$$(n-2)^2+\alpha+1\le (n+1)(n-2)+\alpha-\widehat{\beta}-2\widehat{a}_n-2\le (n-2)^2+\alpha+3$$
and since $\delta\ge n\ge 6$ it follows that either $\gamma=\alpha$, or $\gamma=\alpha-3$ work.
It remains the case  $\widehat{\beta}=0$ where $\beta=0$ or $\beta=n$.
By substituting in the hypothesis of our lemma
the values $(a_{n+1},a_n,a_{n-1})=(0,n-4,2)$ we get
$\beta=n$ and
$\alpha=1$. In this case
 the maximal degree of the trace $X_M\cap M$ is
$$(n+1)(n-2)+1-2(n-3)-2=(n-2)^2-1+n$$
and since $\delta\ge n+1$
we can choose $\gamma=\alpha=1$.

In the last case $(a_n,a_{n-1})=(n-5,3)$, if $n-3\le \widehat{\beta}$, then
 the maximal degree of the trace $X_M\cap M$ is
$$(n-2)^2+\alpha\le (n+1)(n-2)+\alpha-\widehat{\beta}-2\widehat{a}_n-3\le (n-2)^2+\alpha+4$$
and since $\delta\ge n\ge 7$ it follows that either $\gamma=\alpha$, or $\gamma=\alpha-3$ work.
It remains the case  $\widehat{\beta}=0$ where $\beta=0$ or $\beta=n$.
By substituting in the hypothesis of our lemma
the values $(a_{n+1},a_n,a_{n-1})=(0,n-5,3)$ we get
$\beta=n$ and
$\alpha=0$, which is a contradiction.
Then this case is impossible.

Now assume that $a_{n+1}\neq 0$.
In this case we have also $\beta=0$, hence it follows $\widehat{\beta}=0$ and
$\widehat{a}_j=a_j$ for any $4\le j\le n+1$.
By assumption we have
$$3a_{n+1}+2a_{n}+a_{n-1}\ge 3n-7$$
and, as in the first case, we also have
$$a_{n+1}+a_n+a_{n-1}\le n-2$$

These two inequalities imply that $(a_{n+1},a_n,a_{n-1})$ lies in the tetrahedron
with vertices
$(n-2,0,0)$, $(n-3,1,0)$, $(n-\frac{7}{3},0,0)$, $(n-\frac{5}{2}, 0, \frac{1}{2})$.
The only integer points in this tetrahedron are $(n-2,0,0)$ and $(n-3,1,0)$.

In the case $(n-2,0,0)$ the maximal degree of the trace $X_M\cap M$ is
$$(n+1)(n-2)+\alpha-3(n-2)=(n-2)^2+\alpha$$
and clearly the minimal degree is $(n-2)^2$, thus
one of the choices $\gamma=\alpha$ or $\gamma=\alpha-3$ works.
In the case $(n-3,1,0)$  the maximal degree of the trace $X_M\cap M$ is
$$(n+1)(n-2)+\alpha-3(n-3)-2=(n-2)^2+\alpha+1$$
and the minimal degree is obviously $(n-2)^2$,
so that
one of the choices $\gamma=\alpha$ or $\gamma=\alpha-3$ works.
\qedd

\begin{lemma0} \label{terzo}
If all the assumptions of \lemref{primo}
and \lemref{secondo} are not satisfied, then
there exists $\gamma'\ge 0$ satisfying $\gamma'+2\le n-4$,
and every
$\gamma\in[\gamma',\gamma'+2]$ can be attained
by a convenient specialization of type $(a)$.
\end{lemma0}

{\it Proof.} The maximal degree of the trace $X_M\cap M$ is
$${\mathcal M}:=(n+1)(n-2)+\alpha-\beta-3a_{n+1}-2a_n-a_{n-1}$$ .

Since the assumption of Lemma \ref{secondo} are not satisfied,
we have ${\mathcal M}\ge (n-2)^2+2$.

The minimal possible degree of the trace $X_M\cap M$ is
$${m}:=\sum_{i=4}^{n+1}(i-3) a_i + \min\{1,h\}=
(n+1)(n-2)+\alpha-\beta-3\sum_{i=4}^{n+1} a_i+\min\{1-h,0\}\le$$
$$\le (n+1)(n-2)-3\sum_{i=4}^{n+1} a_i \le
(n+1)(n-2)-3(n-2)= (n-2)^2
$$
where we use the fact that $\sum_{i=4}^{n+1} a_i \ge n-2$.
This is true because either $a_{n+1}= n-2$
or $a_{n+1}\le n-3$ and we have
$$\sum_{i=4}^{n} a_i\ge \frac{(n+1)(n-2-a_{n+1})+\alpha-\beta}{n}> n-2-a_{n+1}-1.$$

Hence if ${\mathcal M}\le n-4$ we choose  $\gamma'={\mathcal M}-(n-2)^2-2$.
Otherwise if ${\mathcal M}\ge n-3$  we choose $\gamma'= n-6$.

Both cases work because of the assumption
$${\mathcal M}-{m}=a_n+2a_{n-1}+3\sum_{i=4}^{n-2} a_i-\min\{1-h,0\}\ge 2.$$
\qedd

We can now prove the last preliminary proposition.
Recall that we denote by $X_L^i$ for $i=1,2,3$  the union of the irreducible components $\xi$ of
$X_L$ such that $\deg(\xi:L)=i$.

\begin{prop0}\label{proprep3}
Let $n\ge 5$ and let ${L}\subset\P n$ be  a subspace of codimension
three. Let $X=X_L\cup X_O$ be a scheme contained in a union of double
points such that $X_L$ is supported on $L$ and is general among the
schemes supported on $L$ and $X_O$ is general. 
Assume that $\deg (X_L:L)+\deg X_O={{n+3}\choose 3}-{{n}\choose 3}= 
\frac{3}{2}n^2+\frac{3}{2}n+1$, and that $\deg X_O= (n+1)^2+\alpha$, 
for $0\le \alpha \le n-1$. 
We also assume that the number of the irreducible components of 
$X_L^1\cup X_L^3$ is $\ge \frac{n}{3}$.  
Then there are no cubic hypersurfaces in $\P n$ which contain ${L}$
and which contain $X$. 
\end{prop0}

{\it Proof.}
For $n=5,6,7$ it is a direct computation (see the appendix). 

For $n\ge 8$ the statement follows by induction, and by the sequence
$$0\rig{}I_{{L}\cup{M},{\P n}}(3)\rig{}
I_{{L},{\P n}}(3)\rig{}I_{{L}\cap {M},{M}}(3)\rig{}0$$
where $M$ is a general codimension three subspace. We get
$$0\rig{}I_{X\cup{L}\cup{M},{\P n}}(3)\rig{}
I_{X\cup {L},{\P n}}(3)\rig{}I_{(X\cup{L})\cap {M},{M}}(3).$$

First by Lemmas \ref{primo}, \ref{secondo}, \ref{terzo} we can specialize $X_O=X_O'\cup X_M$
in such a way that $\deg X_O'=3(n+1)+\beta$ (we will call in the following $\widehat{\beta}=\beta$),
$X_M$ is supported on $M$ and
$\deg (X_M\cap M)=(n-2)^2+\gamma$,
where $0\le \beta \le 2n-1$, $0\le \gamma\le n-4$,
$\gamma=\alpha \textrm{\ (mod\ }3)$ and
$\alpha -\beta -n\le \gamma \le \alpha$.
Notice also that we have
$\alpha-\beta-\gamma\ge -2n+4$.
It follows that
$$n-2 \le \deg (X_M:M)=3(n-2)+\alpha-\beta-\gamma\le 4n-6$$

Moreover let us specialize $X_L=X_L'\cup X_L''$ where $X_L'$ is
supported on $L\setminus M$ and $X_L''$ is supported on $L\cap M$.
We may also assume that the number of irreducible components of
$(X_L'')^1\cup (X_L'')^3$ is $\ge \frac{n-3}{3}$. We may assume
that
$$2n-5\le \deg (X_L':L)=3(n-2)+\gamma-\alpha\le 3(n-2)$$
 indeed note that
$3(n-2)+\gamma-\alpha=0\textrm{\ (mod\ } 3)$ and there exist at
least $ \frac{n}{3}$ irreducible component in
$(X_L')^1\cup(X_L')^3$. Note that by using the minimal number of
irreducible component in $(X_L')^1\cup(X_L')^3$,  at least $
\frac{n}{3}-1$ components  remain in $X_L''$, preserving our
inductive assumption.
It follows that
$$\deg (X_L':L)+\deg (X_M:M)+\deg X_O'=9(n-1)$$
moreover we have clearly
$$4n-11\le\deg (X_L':L)+\deg (X_M:M)\le6n-12$$
and we may apply \propref{proprep2} and \propref{proprep2bis}, since
the scheme $X_L'\cup X_M \cup X_O'$ satisfies the corresponding
assumptions.
Then we conclude by induction, indeed the scheme
$(X_M\cup X_L'')\cap M$ satisfies
our assumptions with respect to the spaces $M$ and $M\cap L\subset M$.
Precisely we have (by subtraction)
$$\deg ((X_L''\cap M):(L\cap M))+\deg (X_M\cap M)=
\frac{3}{2}(n-3)^2+\frac{3}{2}(n-3)+1,$$
and $\deg (X_M\cap M)= (n-2)^2+\gamma$,
where $0\le \gamma \le n-4$.
\qedd

We are finally in position to give the proof of the main theorem.

\medskip 

{\it Proof of Theorem \ref{cubiche}.}
We fix a codimension three linear subspace $L\subset{\bf P}^n$ and we
prove the statement by induction by using the exact sequence
$$0\rig{}I_{{L},{\P n}}(3)\rig{}I_{\P n}(3)\rig{}I_L(3).$$

We prove the claim by induction on $n$ from $n-3$ to $n$.
By \lemref{casi-iniziali} we know that the theorem holds for $n=2,3,4$.
Let $X$ be a general scheme contained in a collection of double points and
with $\deg X={{n+3}\choose 3}$.

Since $n\geq 5$ we can assume that $X$ contains at most one component of length $\le 3$.
Fix  a codimension three linear subspace $L\subset{\bf P}^n$ and
consider the exact sequence
\begin{equation}
\label{subspace3}
0\rig{}I_{X\cup L}(3)\rig{}I_{X}(3)\rig{}I_{X\cap L,L}(3).
\end{equation}
We want to specialize on $L$ some components of $X$ so that
$\deg (X\cap L)={{n}\choose 3}$ and apply  \propref{proprep3}.

We keep outside $L$ the irreducible components of $X$
starting from the ones with maximal length in such a way that
$\deg X_O= (n+1)^2+\alpha \ge (n+1)^2$ and it is minimal.
We get by construction that $\alpha\le n-1$.
Let $a_i$ be the number of components of $X_L=X\setminus X_O$ of degree $i$ for $i=4,\ldots, n+1$
and let $h$ be the degree of the component of $X$ of length $\le 3$.
Then $\sum_{i=4}^{n+1}i a_i + h ={{n+3}\choose 3}-(n+1)^2-\alpha$.

After the specialization, the minimum degree of the trace $X_L\cap L$ is
$$\sum_{i=4}^{n+1}(i-3) a_i + 1={{n+3}\choose 3}-(n+1)^2-\alpha-h-3\sum_{i=4}^{n+1} a_i+1$$
if $h\ge 1$ or
$$\sum_{i=4}^{n+1}(i-3) a_i={{n+3}\choose 3}-(n+1)^2-\alpha-3\sum_{i=4}^{n+1} a_i$$
if $h=0$.
On the other hand the maximum degree of the trace $X_L\cap L$ is
$${{n+3}\choose 3}-(n+1)^2-\alpha-3a_{n+1}-2a_n-a_{n-1}.$$

We want to prove that ${{n}\choose 3} $ belongs to the range
between the minimum and the maximum of $\deg(X_L\cap L)$.
This is implied by the inequalities
\begin{equation}\label{first}
\alpha+3a_{n+1}+2a_n+a_{n-1}\le\frac{n(n-1)}{2}
\end{equation}
and
\begin{equation}\label{second}
\frac{n(n-1)}{2}\le \alpha+h+3\sum_{i=4}^{n+1} a_i-1,\quad
\mbox{ or }\quad\frac{n(n-1)}{2}\le \alpha+3\sum_{i=4}^{n+1} a_i
\end{equation}
In order to prove the  inequality (\ref{first}),  consider first the case $a_{n+1}\neq 0$. Then $\alpha=0$ and we have
$$a_{n+1}+\frac{2}{3}a_n+\frac{1}{3}a_{n-1}\le\frac{1}{n+1}\sum_{i=4}^{n+1}i a_i=$$
$$=\frac{1}{n+1}\left[{{n+3}\choose 3}-(n+1)^2-h\right]= \frac{n(n-1)}{6}-\frac{h}{n+1}\le \frac{n(n-1)}{6}$$
as we wanted.
If $a_{n+1}=0$ we get
$$2a_n+a_{n-1}+\alpha\le \frac{2}{n}\sum_{i=4}^{n+1}i a_i +\alpha= \frac{2}{n}\left[{{n+3}\choose 3}-(n+1)^2-h-\alpha\right]
+\alpha\le$$
$$\le \frac{2}{n}\left[{{n+3}\choose 3}-(n+1)^2\right]+(n-1)(1-\frac{2}{n})$$
which is $\le \frac{n(n-1)}{2} $ if $n\ge 6$, as we wanted.

In order to prove the inequality (\ref{second}), notice that
$$\sum_{i=4}^{n+1} a_i\ge\frac{1}{n+1}\sum_{i=4}^{n+1} ia_i=\frac{n(n-1)}{6}-\frac{\alpha+h}{n+1}$$
then if $h=0$ we conclude since
$\alpha(1-\frac{3}{n+1})\ge0$,
while if $h\ge 1$ we conclude by the inequality
$(\alpha+h)(1-\frac{3}{n+1})\ge 1$, which is true
 if $\alpha+h\ge 2$, in particular if $\alpha\ge 1$.

Consider the last case $\alpha=0$ and $h\ge 1$.  If $n\neq 2\textrm{\ (mod\ }3)$,
so that $\frac{n(n-1)}{6}$ is an integer, then $X\setminus X_O$ contains at least $\frac{n(n-1)}{6}+1$
irreducible components and this confirms the inequality.  If $n= 2\textrm{\ (mod\ }3)$,
even $\lfloor\frac{n(n-1)}{6}\rfloor$ double points and one component of length $3$ are not enough to cover all
$X\setminus X_O$. Then $X\setminus X_O$ contains at least   $\lfloor\frac{n(n-1)}{6}\rfloor+2$ irreducible components
and again the inequality is confirmed.

Then a suitable specialization of $X_L$ exists such that $\deg (X_L\cap L)= {{n}\choose 3}$.
We denote again by $X_L^i$ for $i=1,2,3$ the union of irreducible components $\xi$ of
$X_L$ such that $\deg(\xi:L)=i$.

In order to apply \propref{proprep3} we need only to show that the irreducible components of
$X_L^1\cup X_L^3$ are at least $\frac{n}{3}$.
If this condition is not satisfied, we show now that it is possible
to choose another suitable specialization such that again $\deg (X_L\cap L)= {{n}\choose 3}$
but the number of irreducible components of
$X_L^1\cup X_L^3$ is $\ge\frac{n}{3}$.
We assume that the number of irreducible components of
$X_L^1\cup X_L^3$ is $\le\frac{n}{3}$.
Indeed we may perform the following operations, that leave the degree of the trace and of the residual both constant.

\begin{itemize}
\item Pull out a component from  $X_L^2$ to $X_L^3$ and push down another component from  $X_L^2$
to $X_L^1$.
\item Pull out a component from  $X_L^2$ to $X_L^3$ and push down a component of $X_L^1$ .
\item Pull out two components from  $X_L^2$ to $X_L^3$ and push down a component from  $X_L^3$
to $X_L^1$.
\end{itemize}

After such operations have been performed, we get that $X_L$ is still a specialization of a subscheme of $X$,
allowing our semicontinuity argument.

If no of the above operations can be performed, then
$X_L^1$ contains only $a_{n-1}$ components of length $n-1$,
$X_L^2$ contains only $a_n'$ components of length $n$
$X_L^3$ contains only $a_n''$ components of length $n$
and $a_{n+1}$ components of length $n+1$.

Then we get
$$\deg (X_L:L)=a_{n-1}+2a_n'+3a_n''+3a_{n+1}=\frac{n(n-1)}{2}-\alpha$$

hence
$$a_{n}'\ge \frac{n(n-1)}{4}-\frac{\alpha}{2}-\frac{3}{2}(a_{n-1}+a_n''+a_{n+1})$$

On the other hand, we have also

$$\deg (X_L\cap L)= {{n}\choose 3}\ge (n-2)( a_{n-1}+a_n'+a_n''+a_{n+1}) \ge $$
$$\ge (n-2)\left[ \frac{n(n-1)}{4}-\frac{\alpha}{2}-\frac{1}{2}(a_{n-1}+a_n''+a_{n+1})\right]  >$$
$$> (n-2)\left[ \frac{n(n-1)}{4}-\frac{n-1}{2}-\frac{n-1}{6}\right]\ge {{n}\choose 3}$$
where the last inequality is true for $n\ge 8$.  This contradiction concludes the proof.

\qedd

\section{Induction}

In order to prove \thmref{main1} we will work by induction on the
dimension and the degree.
In the following lemmas we describe case by case the initial and
special instances, while in Theorem \ref{main} below we present the
general inductive procedure, which involves the differential Horace
method.

\begin{lemma0}\label{n=2}
A general zero-dimensional scheme $X\subset \P{2}$ contained in a
union of double points imposes independent conditions on ${\cal
  O}_{\P{2}}(d)$ for any $d\geq4$, with the only exception of $d=4$
and $X$ given by the collection of $5$ double points.
\end{lemma0}
{\it Proof.}
Assume that $X$ is a general subscheme of a union of double points with
$\deg(X)={{d+2}\choose{2}}$. If $X$ is a collection of double points
the statement follows from the Alexander-Hirschowitz theorem on
$\P{2}$ (for an easy proof see for example \cite[Theorem 2.4]{BO}).

If $X$ is not a collection of double points, fix a hyperplane
$\P{1}\subset\P{2}$. Since $X$ contains at least a component of
length $1$ or $2$, it is clearly always possible to find a
specialization of $X$ such that the trace has degree exactly
$d+1$. Then we conclude by induction from the {\em Castelnuovo
sequence}
$$0\to I_{X:\P{{1}}}(d-1)\to I_X(d) \to I_{X\cap \P{{1}}}(d).$$

Notice that any subscheme of $5$ double points and any scheme
containing $5$ double points impose independent conditions on
quartics, by Remark \ref{perfetto-difettivo}.
\qedd

We give now an easy technical lemma that we need in the following.
\begin{lemma0}\label{metodo}
Assume that $X$ is a general zero-dimensional scheme contained in a
union of double points of $\P{n}$, which contains at least $n-1$
components of length less than or equal to $n$. Then if
$\deg(X)={{n+d}\choose{n}}$ it is possible to specialize some
components of $X$ on a fixed hyperplane $\P{{n-1}}$ in such a way
that $\deg(X\cap\P{{n-1}})={{n-1+d}\choose{n-1}}$.
\end{lemma0}

{\it Proof.}
By assumption there exist at least $n-1$ components
$\{\eta_1,\ldots,\eta_{n-1}\}$ with $\textrm{length}(\eta_i)\leq n$.
Specialize $\eta_1,\dots,\eta_{n-1}$ on the hyperplane $\P{{n-1}}$ in
such a way that the residual of each component is zero.
Then specialize other components so that
$$\delta={{n-1+d}\choose{n-1}}-\deg(X\cap\P{{n-1}})\geq 0$$
is minimal. If $\delta=0$ the claim is proved, so assume
$\delta\geq1$. Obviously we have $\delta< k-1\leq n$, where $k$ is
the minimal length of the components of $X$ which lie outside
$\P{{n-1}}$. Let $\zeta$ be a component with length $k$. Now we
make the first components $\eta_1,\ldots,\eta_{k-1-\delta}$ having
residual $1$ with respect to $\P{{n-1}}$ and we specialize $\zeta$
on $\P{{n-1}}$ with residual $1$. Notice that this is possible
since $0<k-1-\delta\leq n-1$. \qedd

\begin{lemma0}\label{quartic-case}
Fix $3 \leq n \leq 4$.
A general zero-dimensional scheme $X\subset \P{n}$ contained in a
union of double points imposes independent conditions on ${\cal
  O}_{\P{n}}(4)$, with the following exceptions:
\begin{itemize}
\item $n=3$ and either $X$ is the union of $9$ double points, or $X$
  is the union of $8$ double points and a component of length $3$;
\item $n=4$ and $X$ is the union of $14$ double points.
\end{itemize}
\end{lemma0}
{\it Proof.}
If $X$ is a collection of double points, the statement holds by the
Alexander-Hirschowitz theorem.
Assume that $X$ is a scheme with degree
${{n+4}\choose{4}}$ which is not a union of double points.
Let us denote by $D$ the number of double points in $X$ and by $C$
the number of the components with length less than or equal
to $n$.

If $n=3$ and $C=1$, then $D=8$ and $X$ is the exceptional case of the statement.
If $n=3$ and $C=2$, then $D=8$ and the two components $\eta_1$ and
$\eta_2$ with length less than or equal
to $3$  have necessarily length $1$ and $2$. In this case we specialize
$X$ on $\P{2}$ in such a way that the trace is given by the union of
$\eta_1$,  $\eta_2$  and $4$ double points.
Hence we conclude by the Castelnuovo sequence and by induction.
If $C\geq3$, then we denote by $\eta$ the component of $X$ with minimal length.
We specialize $\eta$ on $\P{2}$ in such a way that its residual
is $1$ if $\textrm{length}(\eta)\geq2$, and
$0$ if $\eta$ is a simple point.
Then we apply the construction of Lemma \ref{metodo} on $X\setminus\eta$
(which has at least two components with length less than or equal
to $3$) and we obtain
a trace different from $5$ double points. Hence we conclude by the
Castelnuovo sequence and by induction.

If $n=4$ and $C=2$, then $X$ is given either by the union of $13$
double points, a component of
length $3$ and one of length $2$,
or by the union of $13$ double points, a component of length $4$ and a
simple point. In the first case we specialize $X$ obtaining a trace given by
$8$ double points, a component of length $2$ and a simple point.
Then we conclude by induction as before.
In the second case we cannot use the Castelnuovo sequence
since we would obtain an exceptional case.
In order to conclude we prove that a general union of $13$ double points
and a component of length $4$ imposes independent conditions on quartics.
Indeed we know that there exists a unique quartic hypersurface through
$14$ double points supported at $p_1,\ldots,p_{14}$. This implies
that for any $i=1,\ldots,14$ there is a unique line $r_i$ through
$p_i$ such that $r_1,\ldots,r_{14}$ are contained in a hyperplane.
Then we consider the scheme $Y$ given by the union of $13$ double points
supported at $\{p_1,\ldots,p_{13}\}$ and the component of length $4$
corresponding to a linear space of dimension $3$ which does not
contain $r_{14}$. It is clear that the scheme $Y$ imposes independent
conditions on quartics, then also the scheme given by the union of
$Y$ and a general simple point does the same.

Assume now that $n=4$ and $C=3$. If $D=13$, then we can
degenerate $X$ to one of the previous cases where the components
with length less than or equal to $4$ are two.
If $D=12$, then the remaining three components have length
either $3,3,4$, or $2,4,4$. In these cases we can obtain as a trace $7$
double points and
three components of length either $2,2,3$, or $1,3,3$, and we conclude
by the Castelnuovo sequence.

If $n=4$ and $C\geq4$, we denote by $\eta$ the component of $X$ with
minimal length.
If $\textrm{length}(\eta)=1$ we can degenerate $X$ to a scheme $X'$ where
the components with length less than or equal to $4$ are one less
and we apply the argument to $X'$.
If $2\leq \textrm{length}(\eta)\leq 3$,
then we specialize $\eta$ on $\P{3}$ in such a way that the residual of $\eta$
is $1$.
Then we apply the construction of Lemma \ref{metodo} on $X\setminus\eta$
(which has at least three components with length less than or equal
to $3$) and we obtain a trace different from $8$ double points and a
component of length $3$.
Moreover with this construction we always avoid a residual given by
$7$ double points.
Hence we conclude by the Castelnuovo sequence.
If $\textrm{length}(\eta)=4$, we have only the following possibilities:
$5$ components of length $4$ and $10$ double points,
$10$ components of length $4$ and $6$ double points,
$15$ components of length $4$ and $2$ double points.
In the first two cases we can obtain trace on $\P{3}$
given by $5$ components of length $3$ and $5$ double points, while in
the third case we can obtain a trace equal to $9$ components of length
$3$ and $2$ double points. Then  we conclude by the Castelnuovo sequence.
\qedd

\begin{lemma0}\label{quartic-case-altri}
Fix $5 \leq n \leq 9$.
A general zero-dimensional scheme $X\subset \P{n}$ contained in a
union of double points imposes independent conditions on ${\cal
  O}_{\P{n}}(4)$.
\end{lemma0}

{\it Proof.}
If $X$ is a collection of double points, the statement holds by the
Alexander-Hirschowitz theorem.
Assume that $X$ is a scheme with degree
${{n+4}\choose{4}}$ which is not a union of double points.
Let us denote by $D$ the number of double points in $X$ and by $C$
the number of the components with length less than or equal
to $n$.

If $n\in\{5,6,8\}$ and $C=2$, then we conclude by degenerating $X$ to a union
of double points.

If $n=5$ and $C=3$, then we get either $D=20$, or $D=19$.
In the first case we conclude degenerating $X$ to the
union of $21$ double points.
In the second case the remaining three components have length $2,5,5$,
or $3,4,5$, or $4,4,4$. Then we can obtain a trace equal to $12$
double points and
three components of length respectively $2,4,4$ in the first case,
or $3,3,4$ in the second and third cases.
Then we conclude by induction.

If $n=5$ and $C=4$, then we have $D\in\{20,19,18\}$.
In the first case we can degenerate $X$ to a union of $21$ double points.
If $X$ can be degenerate to a scheme which contains only three components
with length less than or equal to $5$, we conclude by using the previous results.
Then we have to consider only the two cases where $X$ is given by
$18$ double points and four components of length either $3,5,5,5$,
or $4,4,5,5$.
In these cases we can obtain a trace equal to $12$ double points and
three components of length respectively $2,4,4$ in the first case,
and $3,3,4$ in the second case. Hence we conclude by induction.

If $n=5$ and $C\geq5$, we denote by $\eta$ the component with
minimal length.
Then we specialize $\eta$ on $\P{4}$ in such a way that the residual of $\eta$
is $1$ if $\eta$ if $\textrm{length}(\eta)\geq2$, and
$0$ if $\eta$ is a simple point.
Then we apply the construction of Lemma \ref{metodo} on $X\setminus\eta$
(which has at least four components with length less than or equal
to $5$) and we obtain a trace different from $14$ double points. Hence
we conclude by the Castelnuovo sequence and by induction.

If $n=6$ and $D\geq 21$, we specialize $21$ double points on $\P{5}$
and we conclude by the Castelnuovo sequence.
If $D<21$, then we have $C\geq 5$ and we can apply Lemma \ref{metodo},
concluding  by the Castelnuovo sequence.

If $n=7$ and $D\geq 30$, we specialize $30$  double points on $\P{6}$
and we conclude by the Castelnuovo sequence.
If $D<30$, then we have $C\geq 6$ and we can apply Lemma \ref{metodo}.

If $n=8$ and $C=3$, then either $D=58$ and $X$ can be degenerated to
the union of $59$ double points,
or  $D=57$. In this case the remaining three components can have
length $5,5,8,$ or $5,6,7$, or $6,6,6$. In all these case we can
obtain a trace on $\P{7}$ given by $40$ double points and two
components of total degree $10$.

If $n=8$ and $C=4$ and $X$ can be degenerated to a scheme with less than
$4$ components with length less than or equal to $8$, then we conclude.
Then we have only to consider the case where $D=56$ and the remaining four
components of $X$ have length $3,8,8,8$, or $4,7,8,8$, or $5,6,8,8$,
or $5,7,7,8$, or $6,6,7,8$, or $6,7,7,7$.
In all these cases we obtain a trace on $\P{7}$ given by $40$ double points and
two components of total degree $10$, with the exception of the last case,
where we can obtain a trace given by $39$ double points and three
components of total degree $18$.

If $n=8$ and $C=5$ and $X$ can be degenerated to a scheme with less than
$5$ components with length less than or equal to $8$, then we conclude.
Hence we have only to consider the cases $D=56$ or $D=55$.
Listing all the possible lengths of the remaining five components we
easily notice that we can always obtain a trace on $\P{7}$
given either by $40$ double points and two components of total degree $10$,
or by $39$ double points and three components of total degree $18$.

If $n=8$ and $C=6$ and $X$ can be degenerated to a scheme with less than
$6$ components with length less than or equal to $8$, then we conclude.
Hence we have only to consider the cases $D=55$ or $D=54$.
Listing all the possible lengths of the remaining six components, we
easily notice, as before, that we can always obtain a trace on $\P{7}$
given either by $40$ double points and two components of total degree $10$,
or by $39$ double points and three components of total degree $18$.

If $n=8$ and $C\geq 7$, we apply Lemma \ref{metodo} and we conclude by the
Castelnuovo sequence.

If $n=9$ and $D\geq59$, we specialize $59$ double points on $\P{8}$
and we conclude by the Castelnuovo sequence.
If $D<59$, then we get $C\geq 8$ and we conclude by applying Lemma
\ref{metodo} and by the Castelnuovo sequence.
\qedd

\begin{lemma0}\label{quintic-sestic-case}
Fix $3 \leq n \leq 4$ and $5 \leq d \leq 6$.
A general zero-dimensional scheme $X\subset \P{n}$ contained in a
union of double points imposes independent conditions on ${\cal
  O}_{\P{n}}(d)$.
\end{lemma0}
{\it Proof.}
If $X$ is a collection of double points, the statement holds by the
Alexander-Hirschowitz theorem.
Assume that $X$ is a scheme with degree
${{n+d}\choose{n}}$ which is not a union of double points.

If $(n,d)\neq(4,5)$ and $X$ has only $2$ components with length
less than or equal to $n$, we conclude by degenerating $X$ to a union
of double points.

If $(n,d)=(3,5)$ and $X$ contains at least $7$ double points, we
specialize them on the trace and
we conclude by the Castelnuovo sequence, since the residual contains
$7$ simple points.
If $X$ has less than $7$ double points, then $X$ has obviously at
least $3$ components with length
less than or equal to $3$. In this case we specialize a component with
minimal length
making it having residual $1$,
then we apply the construction of Lemma \ref{metodo} on the remaining
components
and we conclude by the Castelnuovo sequence, since the residual contains
at least a simple point.

If $(n,d)=(4,5)$ and $X$ contains at least $14$ double points, we
specialize them on the trace and
we conclude by the Castelnuovo sequence, since the residual contains
$14$ simple points.
If $X$ has less than $14$ double points, then $X$ has obviously at
least $4$ components with length
less than or equal to $4$. In this case we specialize a component with
minimal length
making it having residual $1$,
then we apply the construction of Lemma \ref{metodo} on the remaining
components
and we conclude by the Castelnuovo sequence, since the residual contains
at least a simple point.

If either $(n,d)=(3,6)$, or $(n,d)=(4,6)$  and $X$ has at least $3$
components with length
less than or equal to $3$, we conclude by Lemma \ref{metodo} and by induction.
\qedd

We are now in position to give the general inductive argument which
completes the proof of Theorem \ref{main1}.

Given a scheme $X\subseteq \P n$ of type $(m_1,\ldots ,m_{n+1})$ and a
fixed hyperplane $\P{{n-1}}\subseteq\P{n}$, we denote for any $1\leq
i\leq n+1$:
\begin{itemize}
\item by $m_i^{(1)}$ the number of component of length $i$ completely
  contained in $\P{{n-1}}$,
\item by $m_i^{(2)}$ the number of
  component of length $i$ supported on $\P{{n-1}}$ and with residual
  $1$ with respect to $\P{{n-1}}$, and
\item by $m_i^{(3)}$ the number
  of component of length $i$ whose support does not lie in
  $\P{{n-1}}$.
\end{itemize}
Obviously we have $m_i^{(1)}+m_i^{(2)}+m_i^{(3)}=m_i$, and
$m_{n+1}^{(1)}=0$, $m_{1}^{(2)}=0$. We denote
$t_i=m_i^{(1)}+m_{i+1}^{(2)}$, for $i=1,\ldots,n+1$,
$r_1=m_1^{(3)}+\sum m_i^{(2)}$, and $r_i=m_i^{(3)}$ for
$i=2,\ldots,n+1$.

\begin{thm0}\label{main}
Fix the integers $n\geq 2$ and $d\geq4$.
A general zero-dimensional scheme $X\subset \P{n}$ contained in a  union of double points
imposes independent conditions on ${\cal O}_{\P{n}}(d)$ with the following
exceptions
\begin{itemize}
\item $n=2$, $d=4$ and $X$ is the union of $5$ double points;
\item $n=3$ and either $X$ is the union of $9$ double points, or $X$
  is the union of $8$ double points and a component of length $3$;
\item $n=4$ and $X$ is the union of $14$ double points.
\end{itemize}
\end{thm0}

{\it Proof.}
We prove the statement by induction on $n$ and $d$.
In Lemma \ref{n=2} we have proved the statement for $n=2,d\geq4$,
in Lemma \ref{quartic-case} and Lemma \ref{quartic-case-altri} for $d=4, 3\leq n\leq 9$
and in Lemma \ref{quintic-sestic-case} for $d=5, n=3,4$ and
$d=6,n=3,4$. Then we need to prove the remaining cases.
Assume $n\geq 3$ and in particular when $d=4$ assume $n\geq 10$, and when $5\leq d\leq 6$
assume $n\geq5$.

The proof by induction is structured as follows:
\begin{itemize}
\item
for $d=4$ and $n\geq 10$, we assume that any scheme in $\P{n}$ imposes
independent conditions on ${\cal O}_{\P{{n-1}}}(4)$.
Recall that any scheme in $\P{n}$ imposes independent conditions on
${\cal O}_{\P{n}}(3)$ (by Theorem \ref{cubiche}) and any scheme of
degree greater than or equal to $(n+1)^2$ imposes independent
conditions on ${\cal O}_{\P{n}}(2)$ (by Theorem \ref{quadrics}).
Then we prove the statement for $d=4, n \geq 10 $;

\item
for $d\geq5$ we assume that any scheme in $\P{a}$ imposes independent conditions on
${\cal O}_{\P{a}}(b)$ for $(a,b)\in\{(n-1,d),(n,d-1),(n,d-2)\}$
and we prove it for $(a,b)=(n,d)$.
\end{itemize}

It is enough to prove the statement for a scheme $X$ with degree
$\deg X={{d+n}\choose{n}}$.

Let $X\subseteq {\bf P}^n$ be a scheme of type
$(m_1,\ldots,m_{n+1})$ contained in a union of double points and
suppose $\deg X=\sum i m_i={{d+n}\choose{n}}$.
Fix a hyperplane $\P{{n-1}}$ in $\P{n}$.

If there is a choice $m_i^{(1)},m_i^{(2)},m_i^{(3)}$ (and $t_i,r_i$ as above) such that
$\sum i t_i={{d+n-1}\choose{n-1}}$, we can conclude by the Castelnuovo sequence and
by induction
$$0\to I_{X:\P{{n-1}}}(d-1)\to I_X(d) \to I_{X\cap \P{{n-1}}}(d).$$
Suppose now that such a decomposition is impossible, and fix a
decomposition $m_i^{(1)},m_i^{(2)},m_i^{(3)}$ such that
${{d+n-1}\choose{n-1}}- \sum i t_i > 0$
is minimal and define
\begin{equation}\label{def-epsilon}
\varepsilon:={{d+n-1}\choose{n-1}}- \sum i t_i.
\end{equation}
Obviously $0< \varepsilon<n$ and $\varepsilon < \min\{i-1: m_i^{(3)}\neq0\}$.
By the minimality of our decomposition we have $m_1^{(3)}=m_2^{(3)}=0$.
Notice that by minimality of $\varepsilon>0$ we have also $m_i^{(2)}=0$ for all $i\neq n+1$.

Now let us define
$$\varepsilon_{n+1}=\min\{\varepsilon,m_{n+1}^{(3)}\},\quad
\varepsilon_{n}=\min\{\varepsilon-\varepsilon_{n+1},m_{n}^{(3)}\}$$
and, for any $i=n-1,\ldots,1$,
$$\varepsilon_{i}=\min\{\varepsilon-\sum_{k=i+1}^{n+1}\varepsilon_{k},m_{i}^{(3)}\}.$$
Obviously we have $\varepsilon_{1}=\varepsilon_{2}=0$
and $\sum_{i=3}^{n+1} \varepsilon_i=\varepsilon$.

{\bf Step 1:} Let $\Gamma\subseteq {\bf P}^{n-1}$ be a general scheme of type
$(\varepsilon_2,\varepsilon_3,\ldots,\varepsilon_{n+1},0)$ supported on a collection
$\{\gamma_1,\ldots,\gamma_{\varepsilon}\}\subseteq{\bf P}^{n}$ of
points and $\Sigma\subseteq{\bf P}^{n}$ a general scheme of type
$(0,0,m_{3}^{(3)}-\varepsilon_{3},\ldots,m_{n+1}^{(3)}-\varepsilon_{n+1})$
supported at points which are not contained in ${\bf P}^{n-1}$.

By induction we know that
$$h_{{\bf P}^{n}}(\Gamma\cup\Sigma,d-1)=\min\left(\deg\Gamma\cup\Sigma,{{n+d-1}\choose{n}}\right)$$
where $\deg \Gamma\cup\Sigma=\sum(i-1)\varepsilon_i+\sum i
(m_i^{(3)}-\varepsilon_i)=\sum im_i^{(3)}-\varepsilon$.

{}From the definition of $\varepsilon$ it follows that
${{n+d-1}\choose{n}}=m_{n+1}^{(2)}+\sum im_i^{(3)}-\varepsilon$
and since of course $m_{n+1}^{(2)}\geq0$, we obtain $h_{{\bf
P}^{n}}(\Gamma\cup\Sigma,d-1)=\sum im_i^{(3)}-\varepsilon$

{\bf Step 2:}
Now we want to add
a collection $\Phi$ of $m_{n+1}^{(2)}$ simple points in $\P{{n-1}}$
to the scheme $\Gamma\cup\Sigma$ and we want to obtain a $(d-1)$-independent scheme.
From the previous step it is clear that
$\dim I_{\Gamma\cup\Sigma}(d-1)=m_{n+1}^{(2)}$.
Hence we have only to prove that there exist no hypersurfaces of degree $d-2$ through $\Sigma$.
Notice that for $d\geq5$ we have
\begin{equation}\label{prima}
\deg(\Sigma)=\sum i (m_i^{(3)}-\varepsilon_i) \geq {{n+d-2}\choose{n}}
\end{equation}
and for $d=4$ and $n\geq 10$ we have
\begin{equation}\label{seconda}
\deg(\Sigma)=\sum i (m_i^{(3)}-\varepsilon_i) \geq (n+1)^2\geq {{n+2}\choose{n}}.
\end{equation}
Indeed by definition of $\varepsilon$, we have
$$\sum i (m_i^{(3)}-\varepsilon_i)={{n+d-1}\choose{n}}+\varepsilon-\sum i\varepsilon_i-m_{n+1}^{(2)}$$
and since
$$\sum i\varepsilon_i-\varepsilon=\sum (i-1)\varepsilon_i\leq n\varepsilon\leq(n-1)n
\quad\textrm{ and }\quad
m_{n+1}^{(2)}\leq\frac{1}{n}{{n+d-1}\choose{n-1}}$$
we obtain
$$\sum i (m_i^{(3)}-\varepsilon_i)\geq
{{n+d-1}\choose{n}}-(n-1)n
-\frac{1}{n}{{n+d-1}\choose{n-1}}=:S(n,d).$$
It is easy to check that for any $d\geq5$ and $n\geq3$ we have $S(n,d)>{{n+d-2}\choose{n}}$,
which proves inequality \eqref{prima}. On the other hand one can also check that
$S(n,4)>(n+1)^2$ for any $n\geq 10$, proving thus inequality \eqref{seconda}.

Then by induction we know that $\Sigma$ imposes independent conditions on ${\cal O}_{\P{n}}(d-2)$,
and so we get $\dim I_{\Sigma}(d-2)=0$.
Thus we obtain
$$h_{{\bf P}^{n}}(\Gamma\cup\Sigma\cup\Phi,d-1)=\sum im_i^{(3)}-\varepsilon+
m_{n+1}^{(2)}={{n+d-1}\choose{n}}.$$

{\bf Step 3:} Let us choose a family of general points
$\{\delta^1_{t_1},\dots,\delta_{t_\varepsilon}^{\varepsilon}\}\subseteq{\bf
  P}^n$, with parameters  $(t_1,\ldots,t_\varepsilon)\in K^\varepsilon$, such that
  for any $i=1,\ldots,\varepsilon$ we have $\delta_0^i=\gamma_i\in{\bf P}^{n-1}$ and
$\delta_{t_i}^i\not\in{\bf P}^{n-1}$  for any $t_i\neq0$.

 Now let us consider the following schemes:

$\bullet\quad\Delta_{(t_1,\ldots,t_\varepsilon)}$ a family of schemes of type
$(\varepsilon_2,\ldots,\varepsilon_{n+1},0)$ supported at the points
$\{\delta^1_{t_1},\dots,\delta_{t_{\varepsilon}}^{\varepsilon}\}$;

$\bullet\quad {\Phi}^2$ of type
$(0,\ldots,0,m_{n+1}^{(2)})$, the union of double points
supported on $\Phi\subset{\bf P}^{n-1}$;

$\bullet\quad\Psi\subseteq{\bf P}^{n-1}$ of type
$(m_1^{(1)},\ldots,m_{n}^{(1)},0)$ supported at general points of
${\bf P}^{n-1}$;

$\bullet\quad$ $\Sigma \subseteq{\bf P}^{n}$, defined in Step 1,
of type $(0,0,m_{3}^{(3)}-\varepsilon_{3},\ldots,m_{n+1}^{(3)}-\varepsilon_{n+1})$.

\smallskip
By induction the scheme $(\Psi\cup \Phi^2|_{\P{{n-1}}}
\cup\Gamma)\subseteq{\bf P}^{n-1}$ has Hilbert function
$$h_{{\bf P}^{n-1}}(\Psi\cup \Phi^2|_{\P{{n-1}}}\cup\Gamma,d)=
\sum i m_{i}^{(1)}+n m_{n+1}^{(2)}+\varepsilon=
\sum i t_{i}+\varepsilon=
{{d+n-1}\choose{n-1}}$$ i.e.\ it is $d$-independent.

In order to prove that $X$ imposes independent conditions on $\mathcal{O}_{\mathbf{P}^n}(d)$,
it is enough to prove the following claim.

{\bf Claim:} There exist $(t_1,\ldots,t_\varepsilon)$ such that the
scheme $\Delta_{(t_1,\ldots,t_\varepsilon)}$ is independent with
respect to the system $I_{\Psi\cup\Phi^2\cup\Sigma}(d)$.

Assume by contradiction that the claim is false. Then by Lemma
\ref{curvilinear} for any $(t_1,\ldots,t_\varepsilon)$ there exist pairs
$(\delta_{t_i}^i,\eta_{t_i}^i)$ for all $i=1,\ldots,\varepsilon$, with $\eta_{t_i}^i$
a curvilinear scheme supported at $\delta_{t_i}^i$ and contained in
$\Delta_{(t_1,\ldots,t_\varepsilon)}$ such that
\begin{equation}\label{assurdo}
h_{{\bf P}^n}(\Psi\cup\Phi^2\cup\Sigma\cup\eta_{t_1}^1\cup\ldots,\eta_{t_\varepsilon}^{\varepsilon},d)<
{{d+n}\choose{n}}-\sum (i-2) \varepsilon_i.
\end{equation}
Let $\eta_0^i$ be the limit of $\eta_{t_i}^i$, for
$i=1,\ldots,\varepsilon$.

Suppose that $\eta_0^i\not\subset{\bf P}^{n-1}$ for $i\in
F\subseteq\{1,\ldots,\varepsilon\}$ and $\eta_0^i\subset{\bf
P}^{n-1}$ for $i\in G=\{1,\ldots,\varepsilon\}\setminus F$.

Given $t\in K$, let us denote $Z_t^F=\cup_{i\in F}(\eta^i_t)$ and
$Z_t^G=\cup_{i\in G}(\eta^i_t)$. Denote by $\widetilde{\eta_0^i}$
for $i\in F$
the residual of $\eta_0^i$ with respect to ${\bf P}^{n-1}$ and by
$f$ and $g$ the cardinalities respectively of $F$ and $G$.

By the semicontinuity of the Hilbert function and by
(\ref{assurdo}) we get
$$h_{{\bf P}^n}(\Psi\cup\Phi^2\cup\Sigma\cup Z^F_0\cup Z^G_t,d)\leq h_{{\bf P}^n}
(\Psi\cup\Phi^2\cup\Sigma\cup Z^F_t\cup
Z^G_t,d)<{{d+n}\choose{n}}-\sum (i-2) \varepsilon_i.$$

On the other hand, by the semicontinuity of the Hilbert function
there exists an open neighborhood $O$ of 0 such that for any $t\in
O$
$$h_{{\bf P}^n}(\Phi\cup\Sigma\cup(\cup_{i\in F}\widetilde{\eta_0^i})\cup Z^G_t,d-1)
\geq h_{{\bf P}^n}(\Phi\cup\Sigma\cup(\cup_{i\in F}\widetilde{\eta_0^i})\cup Z^G_0,d-1)$$
Since the scheme $\Phi\cup\Sigma\cup(\cup_{i\in F}\widetilde{\eta_0^i})\cup Z^G_0$ is contained in
$\Phi\cup\Sigma\cup\Gamma$, which is $(d-1)$-independent by Step 2, we have
$$h_{{\bf P}^n}(\Phi\cup\Sigma\cup(\cup_{i\in F}\widetilde{\eta_0^i})\cup Z^G_0,d-1)
=m_{n+1}^{(2)}+\sum i (m_i^{(3)}-\varepsilon_i)+f+2g.$$

Since
$\Psi\cup \Phi^2|_{\P{{n-1}}}\cup(\cup_{i\in F}\gamma_i)$
is a subscheme of
$\Psi\cup \Phi^2|_{\P{{n-1}}}\cup\Gamma$, which is $d$-independent by Step 3,
it follows that
$$h_{{\bf P}^{n-1}}(\Psi\cup \Phi^2|_{\P{{n-1}}} \cup
(\cup_{i\in F}\gamma_i),d)=\sum i m_i^{(1)}+n m_{n+1}^{(2)}+f$$

Hence for any $t\in O$, by applying the Castelnuovo exact sequence
to the scheme $\Psi\cup \Phi^2 \cup\Sigma\cup Z^F_0\cup
Z^G_t$, we get
$$h_{{\bf P}^n}(\Psi\cup \Phi^2\cup\Sigma\cup Z^F_0\cup Z^G_t,d)\geq$$
$$\geq h_{{\bf P}^n}(\Phi\cup\Sigma\cup(\cup_{i\in F}\widetilde{\eta_0^i})\cup Z^G_t,d-1)+
h_{{\bf P}^{n-1}}(\Psi\cup \Phi^2|_{\P{{n-1}}}\cup
(\cup_{i\in F}\gamma_i),d)\geq$$
$$\geq ( m_{n+1}^{(2)}+\sum i (m_i^{(3)}-\varepsilon_i)+f+2g)
+(\sum i m_i^{(1)}+n m_{n+1}^{(2)}+f)=$$
$$=\sum i m_i-\sum i \varepsilon_i + 2\varepsilon={{d+n}\choose{n}}-\sum (i-2)\varepsilon_i,$$
contradicting \eqref{assurdo}.
This completes the proof of the claim.
\qedd

\section{Appendix}

Here we explain how to compute the dimension of the space
$$V_{d,n}(p_1,\ldots ,p_k, A_1,\ldots, A_k)$$
 defined in (\ref{vdn}) in the introduction. 

These computations are performed in  characteristic
$31991$
using the program Macaulay2 \cite{GS},
and consist essentially in checking that several
square matrices, randomly chosen, have maximal rank.
We underline that
if an integer matrix has maximal rank in positive characteristic,
then it has also maximal rank in 
characteristic zero.
Very likely \thmref{main1}  should be true on any infinite field,
but  a finite number of values for the characteristic  (not including $31991$)
require further and tedious checks, that we have not performed.

Assume that $\dim A_i=a_i$ are given and that $\sum_{i=1}^k (a_i+1)={{n+d}\choose n}=\dim R_{d,n}$.
Consider the monomial basis for $R_{d,n}$ as a matrix $T$ of size
${{n+d}\choose n} \times 1$. 
Consider the jacobian matrix $J$ computed at $p_i$,
which has size  ${{n+d}\choose n} \times (n+1)$. 
Choose a random $(n+1)\times a_i$ integer matrix $A$.
We  concatenate $T$ computed at $p_i$ with $J\cdot A$. It results a  matrix of size
${{n+d}\choose n} \times (a_i+1)$.  When $a_i=n$ (this is the case of Alexander-Hirschowitz theorem)
there is no need  to use a random matrix,
and by Euler identity we can simply take the jacobian matrix $J$ computed at $p_i$.
By  repeating this construction for every point, and placing side by side all these matrices,
we get a square matrix of order ${{n+d}\choose n}$. This is the matrix of coefficients
of the system (\ref{stella}), which corresponds to our interpolation problem. 
Then  there is a unique polynomial $f$ satisfying 
 (\ref{stella}) if and only if the above matrix has maximal rank.
We emphasize that this Montecarlo technique provides a proof, and not only a probabilistic proof.
Indeed consider the subset ${\cal S}$ of points $(p_1,\ldots ,p_k, A_1,\ldots, A_k)$ (lying
in a Grassmann bundle, which locally is isomorphic to the product of affine spaces and Grassmannians, hence 
irreducible) such that the corresponding matrix has  maximal rank.
The subset ${\cal S}$ is open and if it is
not empty, because it contains a random point, then it is dense.

In \propref{proprep},  \propref{proprep2}, \propref{proprep2bis},
\propref{proprep3} we need a modification of the above strategy, since the points 
are supported on some given codimension three subspaces.

As a sample we consider the case considered in \propref{proprep2bis}
where
\texttt{l}$=\deg(X_L:L)=10$,
\texttt{m}$=\deg(X_M:M)=14$, and
\texttt{F}$=\deg(X_O)=39$
and we list below the Macaulay2 script.  Given monomial subspaces $L$ and $M$,
we first compute the cubic polynomials
containing $L$ and $M$, founding a basis of $63$ monomials.
Then we compute all the possible partitions of $10$ and $14$ in 
integers from $1$ to $3$ (which are the possible values  of $\deg(\xi :L)$, resp. $\deg(\xi :M)$, 
where $\xi$ is an irreducible component
of $X_L$, resp. $X_M$), and of $39$ in integers from $1$ to $9$
(which are the possible lengths of a subscheme of a double point in ${\bf P}^8$), 
by excluding the cases which can be easily obtained by degeneration.
We collect the results in the matrices \texttt{tripleL}, \texttt{tripleM} and \texttt{XO},
each row corresponds to a partition.
Then for any combination of rows of the three matrices the program computes 
a matrix \texttt{mat} of order $63$ and its rank.
If the rank is different from $63$ the program prints the case. 
Running the script we see that the output is empty, as we want.

\begin{tiny}
\begin{verbatim}
KK=ZZ/31991;
E=KK[e_0..e_8];
--coordinates in P8

f=ideal(e_0..e_8); 
g=ideal(e_0..e_2); 
h=ideal(e_3..e_5); 
T1=f*g*h;
T=gens gb(T1)
--basis for the space of cubics containing 
--L (e_0=e_1=e_2=0) and M (e_3=e_4=e_5=0) 
--T is a (63x1) matrix

J=jacobian(T);
-- J is a (63x9) matrix

--first case: for the other cases of Proposition 4.8 it is enough
--to change to following line
l=10;m=14;F=39; 
 
---start program
tripleL=matrix{{0,0,0}};
for t from 0 to ceiling(l/3) do 
for d from 0 to ceiling(l/2) do 
for u from 0 to 1 do 
     (if (3*t+2*d+u==l) then tripleL=(tripleL||matrix({{t,d,u}})));

tripleM=matrix{{0,0,0}},
for t from 0 to ceiling(m/3) do 
for d from 0 to ceiling(m/2) do 
for u from 0 to 1 do 
     (if (3*t+2*d+u==m) then tripleM=(tripleM||matrix({{t,d,u}})));

XO=matrix{{0,0,0,0,0,0,0,0,0}};
for n from 0 to ceiling(F/9) do 
     (if (9*n+1==F) then XO=(XO||matrix({{n,0,0,0,0,0,0,0,1}})));
(for n from 0 to ceiling(F/9) do
(for o from 0 to ceiling(F/8) do 
     (if (9*n+8*o+2==F) then XO=(XO||matrix({{n,o,0,0,0,0,0,1,0}})))));
(for n from 0 to ceiling(F/9) do
(for o from 0 to ceiling(F/8) do
(for s from 0 to ceiling(F/7) do 
     (if (9*n+8*o+7*s+3==F) then XO=(XO||matrix({{n,o,s,0,0,0,1,0,0}}))))));
(for n from 0 to ceiling(F/9) do
(for o from 0 to ceiling(F/8) do
(for s from 0 to ceiling(F/7) do 
(for e from 0 to ceiling(F/6) do 
(for c from 0 to ceiling(F/5) do 
     (if (9*n+8*o+7*s+6*e+5*c==F) 
          then XO=(XO||matrix({{n,o,s,e,c,0,0,0,0}}))))))));

k=1;
for a from 1 to (numgens(target(tripleL))-1) do
for b from 1 to (numgens(target(tripleM))-1) do
for c from 1 to (numgens(target(XO))-1) do
     (k=k+1,
mat=random(E^1,E^63)*0,
for i from 1 to tripleL_(a,0) do 
(q1=(matrix(E,{{0,0,0}})|random(E^1,E^6)), mat=mat||random(E^3,E^9)*sub(J,q1)),
for i from 1 to tripleL_(a,1) do 
(q1=(matrix(E,{{0,0,0}})|random(E^1,E^6)), mat=mat||random(E^2,E^9)*sub(J,q1)),
for i from 1 to tripleL_(a,2) do 
(q1=(matrix(E,{{0,0,0}})|random(E^1,E^6)), mat=mat||random(E^1,E^9)*sub(J,q1)),
for i from 1 to tripleM_(b,0) do 
(r1=(random(E^1,E^3)|matrix(E,{{0,0,0}})|random(E^1,E^3)),mat=mat||random(E^3,E^9)*sub(J,r1)),
for i from 1 to tripleM_(b,1) do 
(r1=(random(E^1,E^3)|matrix(E,{{0,0,0}})|random(E^1,E^3)),mat=mat||random(E^2,E^9)*sub(J,r1)),
for i from 1 to tripleM_(b,2) do
(r1=(random(E^1,E^3)|matrix(E,{{0,0,0}})|random(E^1,E^3)),mat=mat||random(E^1,E^9)*sub(J,r1)),
for i from 1 to XO_(c,0) do 
(p1=random(E^1,E^9), mat=mat||sub(J,p1)), 
for i from 1 to XO_(c,1) do 
(p1=random(E^1,E^9), mat=mat||sub(T,p1)||random(E^(8-1),E^9)*sub(J,p1)), 
for i from 1 to XO_(c,2) do 
(p1=random(E^1,E^9), mat=mat||sub(T,p1)||random(E^(7-1),E^9)*sub(J,p1)), 
for i from 1 to XO_(c,3) do 
(p1=random(E^1,E^9), mat=mat||sub(T,p1)||random(E^(6-1),E^9)*sub(J,p1)), 
for i from 1 to XO_(c,4) do 
(p1=random(E^1,E^9), mat=mat||sub(T,p1)||random(E^(5-1),E^9)*sub(J,p1)), 
for i from 1 to XO_(c,5) do 
(p1=random(E^1,E^9), mat=mat||sub(T,p1)||random(E^(4-1),E^9)*sub(J,p1)), 
for i from 1 to XO_(c,6) do 
(p1=random(E^1,E^9), mat=mat||sub(T,p1)||random(E^(3-1),E^9)*sub(J,p1)), 
for i from 1 to XO_(c,7) do 
(p1=random(E^1,E^9), mat=mat||sub(T,p1)||random(E^(2-1),E^9)*sub(J,p1)), 
for i from 1 to XO_(c,8) do mat=mat||sub(T,random(E^1,E^9)),
if (rank(mat)!=63)
then (print(tripleL_(a,0),tripleL_(a,1),tripleL_(a,2),tripleM_(b,0),tripleM_(b,1),tripleM_(b,2),
XO_(c,0),XO_(c,1),XO_(c,2),XO_(c,3),XO_(c,4),XO_(c,5),XO_(c,6),XO_(c,7),XO_(c,8))),
          if (mod(k,29)==0) then print(k));
\end{verbatim}
\end{tiny}

All the others scripts are available at the page
\begin{verbatim}
http://web.math.unifi.it/users/brambill/homepage/macaulay.html
\end{verbatim}

\bigskip
{\sf
Maria Chiara Brambilla:\\
Dipartimento di Matematica e Applicazioni per l'Architettura, Universit\`a di Firenze\\
piazza Ghiberti 27,  50122 Firenze, Italy\\
brambilla@math.unifi.it\\

Giorgio Ottaviani:\\
Dipartimento di Matematica U. Dini, Universit\`a di Firenze\\
viale Morgagni 67/A,  50134 Firenze, Italy\\
ottavian@math.unifi.it
}
\end{document}